    \documentclass[12pt]{amsart}
    \usepackage{latexsym}
    \usepackage[all]{xy}
    \usepackage{amsfonts}
    \usepackage{amsthm}
    \usepackage{amsmath}
    \usepackage{amssymb}
\numberwithin{equation}{section}

    \def\sw#1{{\sb{(#1)}}}
    \def\su#1{{\sp{(#1)}}} 
    \def\sv#1{{\sp{\<#1\>}}}
    \def\suc#1{{\sp{(#1)}}}
    \def\zz#1{{\sb{[#1]}}}

    \def\<{{\langle}}
    \def\>{{\rangle}}
    
    \def\eps{\varepsilon}

    \def\note#1{{}}
    \def\can{{\rm can}}
    \def\cocan{\overline{\rm can}}
    \def\note#1{}

    \def\cM{{\mathcal M}}

    \def\cC{{\mathfrak C}}
    
    \def\cD{{\mathfrak D}}
    \def\cE{{\mathfrak E}}

    \def\lhom#1#2#3{{{\rm Hom}\sb{#1-}(#2,#3)}}
    \def\rhom#1#2#3{{{\rm Hom}\sb{-#1}(#2,#3)}}

    \def\Rend#1#2{{{\rm End}\sp{-#1}(#2)}}
    
    \def\Rhom#1#2#3{{{\rm Hom}\sp{-#1}(#2,#3)}}

    \def\Lhom#1#2#3{{{\rm Hom}\sp{#1-}(#2,#3)}}

    \def\can{{\rm can}}

    \def\beq{\begin{equation}}
    \def\eeq{\end{equation}}
    \def\DC{{\Delta_\cC}}
    \def \eC{{\eps_\cC}}
    \def\DD{{\Delta_\cD}}
    \def \eD{{\eps_\cD}}

    \def\DE{{\Delta_\cE}}
    \def \eE{{\eps_\cE}}

    \def\ot{{\otimes}}

    \def\pl{\Pi^{L}}
    \def\pbl{\overline{\Pi}^{L}}
    \def\pr{\Pi^{R}}
    \def\pbr{\overline{\Pi}^{R}}
 
    \def\sba{_{\alpha}}
    \def\spa{^{\alpha}}
    \def\sbb{_{\beta}}
    \def\spb{^{\beta}}
    \def\spab{^{\alpha\beta}}
    \def\spba{^{\beta\alpha}}
    \def\sbab{_{\alpha\beta}}
    
 \def\sumab{\sum_{\alpha,\beta}}
  \def\suma{\sum_{\alpha}}
  \def\sumAB{\sum_{E,F}}
    \def\sumA{\sum_{E}}
    \def\cuc{\varepsilon_{C}}
    \def\cmc{\Delta_{C}}
    \def\re{\psi_{R}}
    \def\le{\psi_{L}}
    \def\rem{\mathcal{M}(\re)_{A}^{C}}
    \def\lem{{_{A}^{C}\mathcal{M}(\le)}}
    \def\Im{\mathrm{Im\,}}
 \def\hphi{\widehat{\varphi}}
 \def\tcan{\widetilde{\can}}
  \def\th{\tilde{h}}
 
    \newcounter{zlist}
    \newenvironment{zlist}{\begin{list}{(\arabic{zlist})}{
    \usecounter{zlist}\leftmargin2.5em\labelwidth2em\labelsep0.5em
    \topsep0.6ex
    \parsep0.3ex plus0.2ex minus0.1ex}}{\end{list}}

    \newcounter{blist}
    \newenvironment{blist}{\begin{list}{(\alph{blist})}{
    \usecounter{blist}\leftmargin2.5em\labelwidth2em\labelsep0.5em
    \topsep0.6ex 
    \parsep0.3ex plus0.2ex minus0.1ex}}{\end{list}}

    \newcounter{rlist}

    \def\Label#1{\label{#1}\ifmmode\llap{[#1] }\else
    \marginpar{\smash{\hbox{\tiny [#1]}}}\fi}
    \def\Label{\label}

   \newtheorem{proposition}{Proposition}[section]
   \newtheorem{lemma}[proposition]{Lemma}
   
   \newtheorem{corollary}[proposition]{Corollary}
   \newtheorem{theorem}[proposition]{Theorem}

   \theoremstyle{definition}
   \newtheorem{definition}[proposition]{Definition}
   \newtheorem{example}[proposition]{Example}

   \theoremstyle{remark}
   \newtheorem{remark}[proposition]{Remark}

   \newcounter{c}
   
   \newcommand{\etyk}[1]{\vspace{-7.4mm}$$\begin{equation}\Label{#1}
   \addtocounter{c}{1}}
   \renewcommand{\]}{\ifnum \value{c}=1 $$\else \end{equation}\fi}
\begin{document}

   \title{The structure of weak coalgebra-Galois extensions }
   \author{Tomasz Brzezi\'nski}
   \author{Ryan B.\ Turner} \author{Adam P.\ Wrightson}
   \address{ Department of Mathematics, University of Wales Swansea,
   Singleton Park, \newline\indent  Swansea SA2 8PP, U.K.}
   \email{ T.Brzezinski@swansea.ac.uk (T Brzezinski)}

   \email{ marbt@swansea.ac.uk (RB Turner)}  \email{ maapw@swansea.ac.uk (AP
Wrightson)}
  \date{November 2004}
   \begin{abstract}
 Weak coalgebra-Galois extensions are studied. A notion of an invertible
 weak entwining structure is introduced. It is proven that, within an invertible
 weak entwining structure, the surjectivity of the canonical map implies
 bijectivity provided the structure coalgebra $C$ is either coseparable or projective
 as a $C$-comodule.
   \end{abstract}
   \maketitle

    \section{Introduction}
 Galois-type extensions of non-commutative algebras play the role of (schemes of)
 non-commutative principal bundles. The study of generalised Galois extensions was
 initiated by Kreimer and Takeuchi in \cite{KreTak:hop} in terms of
  Hopf-Galois extensions
 and beautifully developed and generalised, in particular  by Doi, Takeuchi and Schneider (cf.\ \cite{Sch:pri}). In recent
 years, in view of the role they play in non-commutative geometry, Hopf-Galois 
 extensions went through a series of further generalisations, thus leading to
 the notion of a coalgebra-Galois extension (cf.\ \cite{BrzMaj:coa}, \cite{BrzHaj:coa}),
 and, most recently, a weak coalgebra-Galois extension (cf.\ \cite{Brz:str}, \cite{CaeDeG:mod}), a weak Hopf-Galois extension (see \cite{CaeDeG:Gal} for the definition and 
 \cite{Kad:act} for an action-free characterisation) and a Hopf algebroid Galois extension (cf.\ \cite{Kad:dep}, \cite{Boh:Gal}). It has been realised in
 \cite{Brz:str} that the general algebraic structure
 underlying all these Galois-type extensions is that of a coring, termed a {\em Galois
 coring} (cf.\ \cite{Wis:Gal}). 
 This, in turn, is a special case of a Galois coring without a grouplike
 element or a {\em Galois comodule} introduced in \cite{KaoGom:com} and recently 
 studied in \cite{Brz:gal}, \cite{CaeDeG:com}, \cite{Wis:gal}.
 
 To have a Galois property  means that a certain map, usually called a {\em 
 canonical map}, must be an isomorphism (or, at least, a bijection). One of the most
 useful results in the standard Hopf-Galois theory is the theorem of Schneider that
 states that it is often enough to prove that the canonical map is an epimorphism
 to conclude that it is an isomorphism. This result proves particularly useful in
 constructing explicit examples of Hopf-Galois extensions, of which there has been 
 a plethora recently, in particular within a realm of
  the non-commutative geometry. Schneider's theorem is
 also very natural from the geometric point of view. Hopf-Galois extensions correspond
 to principal bundles. These are given in terms of free actions of Lie groups on 
 manifolds. Freeness means surjectivity of the canonical map. In differential geometry
 the bijectivity then follows by dimension-type arguments that cannot be transferred
 directly in the algebra context. 
 In \cite[Theorem~4.4]{Brz:gal} it has been shown that the Schneider theorem 
 has in fact a coring
 origin. In the (most general, so far) 
 case of Galois comodules it is enough to check that the canonical map
 is a split epimorphism (in a suitable category) to conclude that it is an isomorphism.
 This then has been applied to bijective entwining structures with a coseparable coalgebra to deduce that the surjectivity of the canonical map suffices to prove that 
 there is a coalgebra-Galois extension (see also \cite{SchSch:gen}
 for a different, coring-free, approach). 
 
 The aim of the present paper is to apply \cite[Theorem~4.4]{Brz:gal} to 
 weak-entwining structures of Caenepeel and De Groot \cite{CaeDeG:mod}
 and thus to prove a Schneider-type structure theorem
 for weak coalgebra-Galois extensions (hence, also, weak Hopf-Galois extensions 
 as a special case). The present paper is organised as follows. In Section 2 
 we set up the notation and conventions, give preliminary results on corings 
 and we also give a formulation of 
 \cite[Theorem~4.4]{Brz:gal} over a general commutative ring. This is a
 mild generalisation of  \cite[Theorem~4.4]{Brz:gal} which  was
 originally stated for vector spaces rather than modules, but it can
 be useful for studying most general Galois-type extensions. As the studies of 
 weak entwining structures and weak Hopf algebras require one to have a 
 number of equalities etc., readily available, we devote Section~3 to collecting
 such useful formulae that are needed for calculations in later sections. We also show
 that one can associate a weak coalgebra-Galois extension to any comodule 
 subalgebra of a weak Hopf algebra. In
 Section~4  we determine the proper notion of an invertible weak
 entwining structure (the naive bijectivity of an entwining map forces a weak
 entwining structure to be an entwining structure). Section~5 contains the first main
 theorem: in the case of an invertible weak entwining structure with a coseparable
 coalgebra, suffices it to prove that the canonical map is an epimorphism to prove
 that there is a weak coalgebra-Galois extension. An example includes a weak 
 Hopf-Galois extension by a weak Hopf algebra with bijective antipode. Finally
 in Section~6 we prove the second main result: within an invertible weak
 entwining structure, if a coalgebra is projective as a comodule, then surjectivity
 of the canonical map implies the bijectivity. As a special case one obtains the 
 following weak Hopf algebra generalisation of the Kreimer-Takeuchi theorem 
 \cite[Theorem~1.7]{KreTak:hop}: for a finite dimensional weak Hopf algebra over
 a field, the surjectivity of the canonical map implies its bijectivity.

    \section{Corings and Galois comodules}
    \subsection{Preliminaries on corings and Galois comodules}
    We work over a commutative ring $k$ with a unit. All algebras are over $k$, associative and with a unit. 
    The product in an algebra is denoted by $\mu$ and the unit, both as an element
    and as a map, is denoted by $1$.
    Unadorned tensor product between $k$-modules is over $k$. 
    All coalgebras are over $k$, coassociative and with a counit. In a coalgebra $C$, the coproduct is denoted by $\Delta_C$ and the counit by $\eps_C$.
     For a ring ($k$-algebra) $A$, the category of right $A$-modules and right $A$-linear maps is denoted by $\cM_A$. Symmetric notation is used for left modules.  
     The dual module of a right $A$-module $M$ is denoted by $M^*$, 
     while the dual of a left $A$-module $N$ is denoted by ${}^*N$. 
     The product in the endomorphism ring of a right module (comodule) is given by composition of maps.

Let $A$ be an algebra. A coproduct in an $A$-coring $\cC$ is denoted by  $\DC :\cC\to \cC\ot_A\cC$, and  the counit is denoted by $\eC:\cC\to A$. To
indicate the action of $\DC$ we use the Sweedler sigma notation,
i.e., for all $c\in \cC$,
$$
\DC(c) = \sum c\sw1\ot c\sw 2, \qquad (\DC\ot_A\cC)\circ\DC(c) =
(\cC\ot_A\DC)\circ\DC(c) = \sum c\sw1\ot c\sw 2\ot c\sw 3,
$$
etc. Capital Gothic  letters always denote corings. The category of right $\cC$-comodules and right $\cC$-colinear maps is denoted by $\cM^\cC$. Recall that  any right $\cC$-comodule is also a right $A$-module, and any right $\cC$-comodule map is right $A$-linear. For a right $\cC$-comodule $M$, $\varrho^M:M\to M\ot_A\cC$ denotes a coaction, and $\Rhom \cC MN$ is the $k$-module of $\cC$-colinear maps $M\to N$. On elements $\varrho^M$ is denoted by the Sweedler notation $\varrho^M(m) = \sum m\zz 0\ot m\zz 1$. Symmetric notation is used for left $\cC$-comodules. In particular, the coaction of a left $\cC$-comodule $N$ is denoted by ${}^N\varrho$, and, on elements, by ${}^N\varrho(n) = \sum n\zz{-1}\ot n\zz 0\in \cC\ot_A N$. The same rules of notation apply to comodules over a coalgebra. A detailed account of the theory of corings and comodules can be found in \cite{BrzWis:cor}.

Given right $\cC$-comodules $M$ and  $N$, the $k$-module $\Rhom\cC MN$ is a right module of the endomorphism ring $B=\Rend\cC M$ with the standard action  $fb = f\circ b$, for all $f\in \Rhom \cC MN$, $b\in B$. This defines a functor $\Rhom \cC M- :\cM^\cC\to \cM_B$, which has 
the left adjoint  $-\ot_BM:\cM_B\to \cM^\cC$ (cf.\ \cite[18.21]{BrzWis:cor}).
The counit of the adjunction is given by the evaluation map
$$
\varphi_{N}: \Rhom \cC M N \ot_BM\to N, \quad f\ot m\mapsto f(m),
$$
Similar adjoint functors exist for left 
$\cC$-comodules. In this case the counit is denoted by $\hphi_N$.

  View $\cC$ as a right $\cC$-comodule with the regular coaction $\DC$. 
$M$ is called a {\em Galois (right) 
 comodule} if $M$ is a finitely generated and
projective right $A$-module, and the evaluation
 map  
 $ \varphi_\cC: \Rhom \cC M\cC\ot_BM\to \cC$
is an isomorphism of right $\cC$-comodules.

Equivalently,  Galois comodules are defined as follows.   If $M$ is a finitely generated projective right $A$-module, then $M^*\ot_BM$ is an $A$-coring with the coproduct $\Delta_{M^*\ot_B M}(\xi\ot m) =  \sum_i \xi\ot e^i \ot \xi^i\ot m$, where $\{e^i\in M,\xi^i\in M^*\}$ is a dual basis of $M_A$, and with the counit $\eps_{M^*\ot_B M}(\xi\ot m) =   \xi(m)$ (cf.\ \cite{KaoGom:com}). In view of the isomorphism $\Rhom \cC M\cC \simeq M^* = \rhom A M A$,  $\varphi_\cC$ becomes the {\em canonical} $A$-coring map
$$
\can_M: M^*\ot_B M\to \cC, \qquad \xi\ot m\mapsto \sum \xi(m\zz0)m\zz 1.
$$
$M$ (with $M_A$ finitely generated projective) is a Galois comodule if and only if the canonical map $\can_M$ is an isomorphism of corings. 

A Galois comodule $M$ with the  endomorphism ring $B$ 
is called a {\em principal comodule} if it is projective
as a left $B$-module.

If $M$ is a  right $\cC$-comodule that is finitely generated and projective as a right
$A$-module, then $M^*$ is a left $\cC$-comodule with the coaction  ${}^{M^*}\!\varrho(\xi) = \sum _i \xi(e^i\zz 0)e^i\zz 1\ot \xi^i$, where $\{e^i\in M,\xi^i\in M^*\}$ is a dual basis of $M$. The endomorphism ring of $M^*$ (as a left $\cC$-comodule) is
isomorphic to the endomorphism ring of $M$.
One can develop the theory of left Galois  and principal comodules
along the same lines as for the right comodule case. In particular if $M$ is a right
Galois  comodule, then $M^*$ is a left Galois 
 comodule.

The case in which $A$ is a Galois $\cC$-comodule is of fundamental importance. In this case the coaction $\varrho^A :A\to A\ot_A\cC\simeq \cC$ is fully determined by a group-like element $g = \varrho^A(1)\in\cC$, i.e., $\varrho^A(a) = ga$. The endomorphism ring $B=\Rend \cC A$ coincides with the subalgebra of $g$-coinvariants in $A$, i.e., $B = A^{co\cC}_g := \{b\in A\; |\; bg=gb\}$. Obviously, $A$ is a finitely generated projective right $A$-module, $A^* \simeq A$, and $A\ot_B A$ is the {\em Sweedler $A$-coring}, with coproduct $a\ot a'\mapsto a\ot 1\ot 1\ot a'$ and counit $a\ot a'\mapsto aa'$. The canonical map comes out as
\begin{equation}\label{canA}
\can_A :A\ot_B A\to \cC, \qquad a\ot a'\mapsto aga'.
\end{equation}
Thus $A$ is a Galois comodule if and only if $\cC$ is a {\em Galois coring} with respect to $g$, a notion introduced in \cite{Brz:str}.

    \subsection{Split canonical epimorphisms are isomorphisms}
    
   The aim of the present paper is to apply \cite[Theorem~4.4]{Brz:gal} to the case
   of weak entwining structures and thus to infer that, similarly to the case of invertible
   entwining structures dealt with in  \cite[Theorem~4.6]{Brz:gal}, in many
   cases it is sufficient to 
   prove that the canonical map is surjective to conclude that it is an isomorphism.
   Theorem~4.4 in \cite{Brz:gal} is formulated for algebras over a field as this is the
   most interesting case from the non-commutative geometry point of view, which was
   the main motivation for introducing the notion of a principal comodule. However,
   the methods of proof of Theorem~4.4 can easily be extended to general algebras
   over a commutative ring. For the completeness we give a version of 
   \cite[Theorem~4.4]{Brz:gal} for algebras over a commutative ring.
   \begin{theorem}\label{thm.mod.main}
   Let $\cC$ be an $A$-coring and $M$ a right $\cC$-comodule that is finitely
   generated and projective as a right $A$-module. Let $B=\Rend\cC M$. Suppose
   that:
   \begin{blist}
   \item $\cC$ is a flat left $A$-module;
   \item the map 
   $$
   B\ot M\to \Lhom\cC{M^*}{M^*\ot M}, \qquad b\ot m\mapsto[\xi\mapsto \xi\circ b\ot m]
   $$
   is an isomorphism of left $B$-modules;
   \item the map
   $$
   \widetilde{can}_M:M^*\ot M\to \cC, \qquad \xi\ot m\mapsto \sum \xi(m\zz 0) m\zz 1,
   $$
   is a split epimorphism of left $\cC$-comodules.
   \end{blist}
   Then $M$ is a Galois comodule and $M$ is a $k$-relatively projective left 
   $B$-module (meaning that any $B$-module epimorphism $N\to M$ that splits 
   as a $k$-module map splits as a $B$-module map).
   \end{theorem}
   \begin{proof}
   Note that the assumption (b) incorporates the part of \cite[Lemma~4.5]{Brz:gal}
   that is used for the proof of  \cite[Theorem~4.4]{Brz:gal}. Once this is realised,
   the same method of proof as in \cite[Theorem~4.4]{Brz:gal} can be used. 
   We repeat the main arguments for completeness.
   
   By assumption (c), the coring $\cC$ is a direct summand of $M^*\ot M$ as a left
   $\cC$-comodule, hence, in view of assumption (b), $\Lhom \cC {M^*}\cC \simeq {}^*(M^*) \simeq M$ is a direct summand of a left $B$-module $B\ot M$. Since
   $B\ot M$ is a $k$-relatively projective left $B$-module, so is $M$.

 The counit of the adjunction $\hphi_{M^*\ot M}:M^*\ot_B \Lhom \cC {M^*} {M^*\ot M}\to
 M^*\ot M$ factorises through the isomorphism in assumption (b) tensored
 with $M^*$, and through the
 obvious isomorphism $M^*\ot_BB\ot M\to M^*\ot M$. Hence $\hphi_{M^*\ot M}$ is an
 isomorphism.
Next we can consider the following diagram, which is commutative in all possible directions since $\hphi$ is a natural transformation,
$$
\xymatrix{M^*\ot_B \Lhom \cC {M^*} {M^*\ot M}\ar@<1ex>[dd]^{M^*\ot_B \Lhom \cC {M^*} {\tcan_M}}\ar[rrr]^{\hphi_{M^*\ot M}}&&& 
   M^*\ot M\ar@<1ex>[dd]^{\tcan_M}\\
&&&\\
M^*\ot_B \Lhom \cC {M^*} \cC\ar@<1ex>[uu]\ar[rrr]^{\hphi_\cC}&&& 
\cC\ar@<1ex>[uu] .}
$$
The upward pointing arrows are sections of $M^*\ot_B \Lhom \cC {M^*} {\tcan_M}$ and $\tcan_M$ respectively. Since $\hphi_{M^*\ot M}$ is an isomorphism, the map $\hphi_\cC$ is bijective (it is a $k$-linear isomorphism). The identifications $\Lhom \cC {M^*}\cC\simeq M$ and $M^*\simeq \Rhom \cC M\cC$  lead to
a $k$-linear isomorphism $M^*\ot_B \Lhom \cC {M^*}\cC
\simeq  \Rhom \cC M\cC \ot_BM$.
In view of this isomorphism,  
the fact that $\hphi_\cC$ is bijective implies that $\varphi_\cC$ is bijective. By assumption, $\cC$ is a flat left $A$-module, so $\cM^\cC$ is an Abelian category. Since $\varphi_\cC$ is a bijective morphism
 in $\cM^\cC$, it is an isomorphism. Thus $M$ is a Galois right $\cC$-comodule.
 \end{proof}
 
 The assumption (b) in Theorem~\ref{thm.mod.main} is satisfied if $M$ is a flat 
 $k$-module. Thus Theorem~\ref{thm.mod.main} implies the following
  \begin{corollary}\label{cor.mod.main}
   Let $\cC$ be an $A$-coring and $M$ a right $\cC$-comodule that is finitely
   generated and projective as a right $A$-module. Let $B=\Rend\cC M$. Suppose
   that:
   \begin{blist}
   \item $\cC$ is a flat left $A$-module;
   \item $M$ is a projective $k$-module;
      \item the map
   $$
   \widetilde{can}_M:M^*\ot M\to \cC, \qquad \xi\ot m\mapsto \sum \xi(m\zz 0) m\zz 1,
   $$
   is a split epimorphism of left $\cC$-comodules.
   \end{blist}
   Then $M$ is a principal comodule. 
   \end{corollary}
   \begin{proof}
   A projective module is flat, hence assumption (b) in Corollary~\ref{cor.mod.main} implies assumption (b) in
   Theorem~\ref{thm.mod.main}. Since $M$ is a $k$-relatively projective left $B$
   module and $M$ is a projective $k$-module,  $M$  is a projective $B$-module, hence
   a principal comodule as required.
   \end{proof}
   \begin{remark}\label{rem.abelian}
   (1) Note that the assumption that $\cC$ is a flat left $A$-module in 
   Theorem~\ref{thm.mod.main} is made to ensure that $\cM^\cC$ is an
   Abelian category (so that  every bijection in $\cM^\cC$ is an isomorphism). 
   Without this assumption  the arguments
  of the proof of Theorem~\ref{thm.mod.main} imply that $\varphi_\cC$, hence
  also $\can_M$, is a $k$-linear isomorphism. Thus, skipping assumption (a),
  one can prove that $M$ is a {\em Galois module in a weaker sense} considered in
  special cases by some authors (cf.\ \cite{Boh:Gal}, \cite{SchSch:gen}).
  
  (2) In the case of an $A$-coring $\cC$ with a grouplike element $g$ and
  the right $\cC$-comodule $M=A$, the dual $M^*=A$ and
  $$
  \Lhom \cC A {A\ot A} = {}^{co\cC}(A\ot A)_g := \{\sum_i a_i\ot b_i\in A\ot A \; |\; \sum_i ga_i\ot b_i = \sum_i a_ig\ot b_i\}.
  $$
  Hence the map in Theorem~\ref{thm.mod.main}~(a)
  comes out as
  $$
  B\ot A\to  {}^{co\cC}(A\ot A)_g, \qquad b\ot a\mapsto b\ot a.
  $$
  Thus the condition (a) means in this case that
  $$
  B\ot A =  \{\sum_i a_i\ot b_i\in A\ot A \; |\; \sum_i ga_i\ot b_i = \sum_i a_ig\ot b_i\},
  $$
  and appears in this form in \cite{Boh:Gal} and \cite{SchSch:gen}.
   \end{remark}

    \section{Weak entwining structures and weak Hopf algebras}
  \subsection{Weak entwining structures and weak coalgebra Galois extensions}  Motivated by a connection between entwining structures and Doi-Koppinen 
    modules (\cite{Brz:mod}, see also  \cite{CaeMil:gen} and 
    \cite[Chapter~5]{BrzWis:cor} for reviews), weak entwining structures 
    were introduced in \cite{CaeDeG:mod} as a structure behind a weak 
    Doi-Koppinen datum defined in \cite{Boh:Doi}. A {\rm right-right weak 
    entwining structure} is a triple
    $(A,C,\re)$, where  $A$ is a $k$-algebra, $C$ is a $k$-coalgebra, and
    $\re : C \otimes A \to A \otimes C$ is a $k$-linear map, 
which, writing
    $\re(c\otimes a)=\sum\sba a \sba \otimes c \spa$, $\re(c\otimes a)=\sum\sbb a \sbb \otimes c \spb$, etc., satisfies
    the following relations:
    \begin{gather}
    \sum\sba (ab)\sba \otimes c\spa = \sumab a \sba b \sbb \otimes
    c\spab, \label{re1}\\
    \sum\sba a \sba \cuc(c\spa)=\sum\sba \cuc(c\spa)1\sba a, \label{re2}\\
    \sum\sba a \sba \otimes \cmc (c \spa) = \sumab a \sbab \otimes
    c\sw 1 \spb \otimes c\sw 2\spa,  \label{re3}\\
    \sum\sba 1 \sba \otimes c \spa = \sum\sba \cuc(c \sw 1 \spa)1\sba
    \otimes c\sw 2. \label{re4}
    \end{gather}
An example of a right-right weak entwining structure is provided by a 
{\em right-right entwining
structure $(A,C,\psi_R)$}. Recall  from \cite{BrzMaj:coa} 
that this is defined by requiring that 
equations (\ref{re1}) and (\ref{re3}) are satisfied, while 
equations (\ref{re2}) and (\ref{re4}) are replaced by
$$
  \sum\sba a \sba \cuc(c\spa)=\cuc(c)a, \qquad \sum\sba 1 \sba \otimes c \spa = 1
    \otimes c.
  $$

    Similarly a {\em left-left weak entwining structure} over $k$ is a triple
    $(A,C,\le)$ consisting of an algebra $A$, a coalgebra $C$, and a $k$-linear map
    $\le : A \otimes C \to C \otimes A$ such that, writing
    $\le(a\otimes c)=\sumA c _{E} \otimes a ^{E} $, $\le(a\otimes c)=\sum_F c _{F} \otimes a ^{F} $ etc., the following relations
    \begin{gather}
    \sumA c_{E} \otimes (ab)^{E}=\sumAB c_{EF}\otimes a^{F}b^{E},  \label{le1}\\
    \sumA \cuc(c_{E})a^{E} = \sumA a \cuc(c_{E})1^{E},  \label{le2}\\
    \sumA \cmc(c_{E}) \otimes a^{E} = \sumAB c _{(1)E} \otimes c_{(2)F} \otimes
a^{EF},  \label{le3}\\
    \sumA c_{E} \otimes 1^{E}=\sumA c _{(1)}\otimes \cuc
    (c_{(2)E})1^{E}  \label{le4}
    \end{gather}
are satisfied. 

Throughout this paper the repeated lower case Greek indices in a {\em down-up  
Einstein convention}, such as $\suma a_\alpha\ot c^\alpha$, always denote the 
action of a right-right weak entwining structure.  Similarly, the repeated upper-case Latin indices $E,F$ as in $\sumA c_E\ot a^E$ denote the action of a left-left weak entwining structure.

    Let $(A,C,\re)$ be a right-right weak entwining
    structure. A $k$-module $M$ together with a right $A$-action and a
    right $C$-coaction is called a {\em weak entwined module over
    $(A,C,\re)$} if
    \begin{gather}
    \varrho^{M} (ma)=\suma m_{[0]}a \sba \otimes {m_{[1]}}\spa,
    \end{gather}
    for all $m \in M$ and $a \in A$.  The category of weak entwined modules
    and $A$-linear $C$-colinear maps is denoted be $\rem$.
    Similarly, given a left-left weak entwining structure $(A,C,\le)$, a $k$-module $M$ together with a left $A$-action and a
    left $C$-coaction is called a {\em weak entwined module over
    $(A,C,\le)$} if
    \begin{gather}
    ^{M}\varrho (am)=\sumA m_{[-1]E} \otimes a^{E}m_{[0]}
    \end{gather}
    for all $m \in M$ and $a \in A$.  The category of weak entwined modules
    over  $(A,C,\le)$ and $A$-linear $C$-colinear maps is  denoted by $\lem$.

The most natural point of view on weak entwining structures and associated modules
is provided by corings and their comodules\footnote{Note that weak 
entwining structures we discuss in this paper 
are {\em self-dual} entwining structures in terminology of \cite{BrzWis:cor}. 
More general weak entwining  structures discussed in \cite{BrzWis:cor} are best
described in terms of {\em weak corings} \cite{Wis:wea}.}. As shown in \cite{Brz:str}, to any right-right weak entwining structure $(A,C,\psi_R)$ one can associate an  $A$-coring $\cC$ such that the category of right $\cC$-comodules is isomorphic to the category of  weak entwined modules over $(A,C,\psi_R)$.
 As understanding of this fact is crucial to what follows, we quote it in full.
    \begin{proposition}\label{prop.ent.cor}
    Let $(A,C,\re)$ be a right-right weak entwining structure. Let
    $$
    p_{R} : A \otimes C \to A \otimes C, \qquad p_{R}=(\mu \otimes
    C) \circ (A \otimes \re) \circ (A \otimes C \otimes 1),
    $$
    and
    $$
    \cC= \mathrm{Im} \ p_{R}=\{\sum_{i,\alpha} a_{i}1\sba \otimes c_{i}\spa \;
    \vert \; \sum_i a_{i} \otimes c_{i} \in A \otimes C \}
    $$
    Then $p_R$ is a projection, i.e., $p_{R} \circ p_{R} = p_{R}$, and
    \begin{enumerate}
    \item
    $\cC$ is an $(A,A)$-bimodule with the left action
    $a'(\suma a'1\sba \otimes c \spa) = \suma a'a1 \sba \otimes c \spa$ and
    the right action $(\suma a'1\sba \otimes c\spa)  a = \sumab a'1 \sba a
    \sbb \otimes c \spab $.
    \item
    $\cC$ is an $A$-coring with coproduct $\DC=(A \otimes \cmc)
    \vert _{\cC}$ and counit $\eC=(A \otimes \cuc) \vert _{\cC}$. Explicitly, for all 
 $a\in A$, $c\in C$,
 $$
 \DC ( \suma a1\sba \otimes c \spa) = \suma a1\sba \otimes c \spa\sw 1\ot_A 1 
 \ot c^\alpha\sw 2 =
 \sumab a1\sba \otimes c\sw 1 \spa\ot_A 1_\beta \ot c\sw 2^\beta.
 $$
    \item
    $\rem = \mathcal{M}^{\cC}$
    \end{enumerate}
\end{proposition}

In a similar way, to a left-left weak entwining structure
    $(A,C,\le)$, one can associate a projection 
    $$
    p_{L} : C \otimes A \to C \otimes A, \qquad p_{L}=(C \otimes
    \mu) \circ (\le \otimes A) \circ (1 \otimes C \otimes A),
    $$
    and an $A$-coring
    $$
    \cD= \mathrm{Im} \ p_{L} = \{\sum_{i, E} c^i_E \otimes 1^Ea^i \;
    \vert \;\sum_i c^{i} \otimes a^{i} \in C \otimes A \}. 
    $$
In this case $\cD$ is an $(A,A)$-bimodule with the left multiplication $a (\sumA 
c_E\ot 1^E a') = \sumAB c_{EF}\ot a^F 1^E a'$ and the right multiplication 
$(\sumA c_E\ot 1^E a') a = \sumA c_E\ot 1^E a'a$. The coproduct and the 
counit of $\cD$ are obtained by restricting of $\Delta_C\ot A$ and $\eps_C\ot A$,  
respectively, to $\mathrm{Im} \ p_{L}$.

The following example, taken from \cite{Brz:str}, recalls the definition of the main 
object of studies of the present paper.

\begin{example}\label{ex.can} Let $C$ be a coalgebra, $A$ an algebra and a right
$C$-comodule with the coaction $\varrho ^{A}$. Let
$$
B= A ^{\mathrm{co}  C}=\{ b \in A | \ \forall a \in A, \ \varrho^{A}
(ba)=b \varrho ^{A} (a)\},
$$
and let
$$ \cocan : A \otimes _{B} A \to A \otimes C, \qquad a \otimes _{B} a' \mapsto
a\varrho ^{A} (a') $$
View $A \otimes _{B} A$ as a left $A$-module via $\mu \otimes _{B}
A$ and a right $C$-comodule via $A \otimes _{B} \varrho ^{A}$. View
$A \otimes C$ as a left $A$-module via $\mu \otimes C$ and as a
right $C$-comodule via $A \otimes \cmc$.  Now suppose the $\cocan$ is
a split monomorphism in the category of left $A$-modules and right
$C$-comodules, i.e., there exists a left $A$-module, right
$C$-comodule map $\sigma : A \otimes C \to A \otimes _{B} A$ such
that $ \sigma \circ \cocan = A \otimes _{B} A$. Let $\tau : C \to A
\otimes _{B} A$, $c \mapsto \sigma (1 \otimes c)$. Define
$$
\re^\sigma: C \otimes A \to A \otimes C, \qquad \re^\sigma = \cocan \circ (A
\otimes _{B} \mu) \circ (\tau \otimes A)
$$
Then $(A,C,\re^\sigma)$ is a weak entwining structure.  The extension of
algebras $B \subseteq A$ is called a (right) {\em weak coalgebra-Galois $C$-extension} and  $(A,C,\re^\sigma)$ is called a {\em canonical right-right weak entwining structure} associated to the 
weak coalgebra-Galois $C$-extension $B\subseteq A$.
\end{example}

By Proposition~\ref{prop.ent.cor}, any weak coalgebra-Galois $C$-extension
$B\subseteq A$ induces a coring $\cC = \Im p_R$ associated to the canonical 
weak entwining structure $(A,C,\re^\sigma)$. As observed in \cite{Brz:str}, the coring $\cC$ is a Galois 
coring, i.e., $A$ is a Galois comodule. Conversely, suppose $(A,C,\re)$ is a 
right-right entwining structure, and $\cC = \Im p_R$ is the corresponding 
coring. Suppose further that $A$ is a right Galois $\cC$-comodule, and let $B$
be the subalgebra of coinvariants. If $A$ is a Galois right comodule, 
then $B\subseteq A$ 
is a weak coalgebra-Galois $C$-extension and, necessarily, $(A,C,\re)$ is the canonical 
right-right entwining structure associated to this extension. This follows immediately from the commutative diagram of right $A$-module left $C$-comodule maps
$$\xymatrix{
   &A\ot_BA \ar[rr]^\cocan \ar@{=}[d]  &&
     A\ot C\ar[d]^{p_R}  &
    \\
  0\ar[r] & A\ot_BA \ar@<.5ex>[rr]^{\can_A}  &&\cC \ar@<.5ex>[ll]^{\can_A^{-1}}\ar[r]\ar[d] & 0 \\
   & && 0 & 
} $$
where $\can_A$ is the canonical map as in equation~\eqref{canA}.

Thus there are two equivalent points of view on weak coalgebra-Galois extensions: 
the entwining-free definition of Example~\ref{ex.can} or the definition within a 
weak entwining structure (or, more precisely, the associated coring). We will make 
substantial use of this latter point of view.

\subsection{Weak Hopf algebras and weak Hopf-Galois extensions}
The notion of a weak bialgebra was introduced in \cite{BohSzl:coa} and \cite{Nil:axi}. 
The paper \cite{Boh:wea} contains a detailed account of the theory of weak 
bialgebras and weak Hopf algebras and is a gold mine of useful formulae and 
properties of weak bialgebras. A concise review of these properties (for infinite  
weak bialgebras over a commutative ring) can be found in 
\cite[Section~36]{BrzWis:cor}.

A weak $k$-bialgebra $H$ is a $k$-module with a $k$-algebra
structure $(\mu, 1)$ and a $k$-coalgebra structure $(\Delta,
\eps)$ such that $\Delta$ is a multiplicative map and
\begin{gather}
\Delta^{2}(1) = \sum 1 \sw 1 \otimes 1 \sw 2 1 _{(1')} \otimes
1_{(2')}= \sum 1 \sw 1 \otimes 1_{(1')} 1\sw 2 \otimes 1_{(2')}\\
\eps (hkl) = \sum \eps (hk \sw 1) \eps (k\sw 2 l) =
\sum \eps (hk \sw 2) \eps (k\sw 1 l) \label{eps}
\end{gather}
for all $h,k,l \in H$.  Here $\Delta(1) = \sum 1\sw 1\ot 1\sw 2 = 
\sum 1 _{(1')} \otimes 1 _{(2')}$.  An important role in the studies 
of the structure of weak bialgebras is played by 
the projections $\pl$, $\pr$, $\pbl$ and $\pbr$, given by
$$
\pl(h)=\sum \eps(1\sw 1 h)1\sw 2, \qquad \pbl (h)= \sum \eps (1 \sw
2 h) 1 \sw 1, $$
$$\pr (h) = \sum \eps(h 1 \sw 2) 1 \sw 1, \qquad
\pbr (h) = \sum \eps (h 1 \sw 1 ) 1 \sw 2.
$$
Every weak bialgebra can be seen to be a {\em bialgebroid}; the above
 projections are used in the construction of a {\em base algebra}  for this
 bialgebroid. Furthermore, they are needed for the definition of an antipode.

A weak Hopf algebra is a weak $k$-bialgebra with a $k$-linear map
$S : H \to H$, called the {\em antipode}, such that, for all $h\in H$, 
$$ \sum h\sw 1S(h\sw 2) = \pl(h), \quad 
\sum S(h\sw 1) h\sw 2 S(h\sw 3) = S(h), \quad \sum S(h\sw 1)  h\sw 2
=\pr(h). $$
The antipode  $S$ is an anti-algebra and anti-coalgebra map. 
Furthermore, 
\begin{equation}
\pl =\pbr\circ S, \qquad
\pr =\pbl\circ S, \label{pi.s}
\end{equation}
and
\begin{equation}
S\circ \pl = \pr\circ S, \qquad S\circ\pr = \pl\circ S. \label{s.pi}
\end{equation}

Given a weak bialgebra $H$,  a {\em right $H$-comodule algebra} is defined in \cite[Definition~2.1]{Boh:Doi} as a $k$-algebra $A$ with
a right $H$ coaction $\varrho:A\to A\ot H$, such that for all $a,b\in A$,
 \begin{equation}\label{comal.m}
 \varrho (ab)=
\varrho(a) \varrho (b),
 \end{equation}
and
\begin{equation}\label{comal.pl}
\sum a \zz 0 \otimes \pl (a \zz 1) = \sum 1 \zz 0 a \otimes 1 \zz
1.
\end{equation}
As shown in \cite[Definition~2.1]{Boh:Doi} and \cite[Proposition~4.10]{CaeDeG:mod}, 
given the multiplicativity of $\varrho$ (\ref{comal.m}), the condition (\ref{comal.pl}) is 
equivalent to each one
of the following statements (which we list here for future reference):
\begin{gather}
\varrho ^2(1)=\sum 1 \zz 0 \otimes 1 \zz 1 1\sw 1 \otimes 1\sw 2, \label{one.1}\\
\varrho^2(1)=\sum 1 \zz 0 \otimes 1 \sw 1 1 \zz 1 \otimes 1 \sw 2,\label{one.2}\\
\sum a \zz 0 \otimes \pbr (a \zz 1) = \sum a 1 \zz 0 \otimes 1\zz 1,
\label{comal.pr}\\
\sum 1\zz 0 \otimes \pbr (1 \zz 1) = \varrho (1),\label{comal.pr.one}\\
\sum 1\zz 0 \otimes \pl(1 \zz 1) = \varrho (1).\label{comal.pl.one}
\end{gather}

Dually, a {\em right $H$-module coalgebra} is defined as a coalgebra $C$ and 
a right $H$-module such that, for all $h,k\in H$ and $c\in C$,
\begin{equation}\label{modco.com}
\Delta_C(ch) = \sum c\sw 1 h\sw 1\ot c\sw 2 h\sw 2,
\end{equation}
and
\begin{equation}\label{modco.counit}
\eps_C (chk) = \sum \eps_C(ch\sw 2)\eps(h\sw 1k).
\end{equation}
Note that condition (\ref{modco.com})  expresses the comultiplicativity of the action,
and hence it is dual to condition  (\ref{comal.m}). The condition 
(\ref{modco.counit}) is dual to (\ref{one.1}). As is the case for comodule algebras,
 various equivalent formulations of condition (\ref{modco.counit}) are possible 
 (cf.\ \cite[Proposition~4.11]{CaeDeG:mod} for details).
 
Left $H$-comodule algebras and left $H$-module coalgebras can be defined in  a 
similar way.

Comodule algebras and module coalgebras of weak bialgebras provide
one with examples of  weak entwining structures. 
More precisely, given a right $H$-comodule algebra $A$ and a right 
$H$-module coalgebra $C$, one 
defines a $k$-linear map
\begin{equation}\label{doi.ent.r}
\psi_R :C\ot A\to A\ot C, \quad c\ot a \mapsto \sum a\zz 0 \ot ca\zz 1 .
\end{equation}
The triple $(A,C,\psi_R)$ is a right-right weak entwining structure (cf.\ \cite[Theorem~4.14]{CaeDeG:mod}). Similarly, given  a left $H$-comodule algebra 
$A$ and a left
$H$-module coalgebra $C$, the $k$-linear map
\begin{equation}\label{doi.ent.l}
\psi_L :A\ot C\to C\ot A, \quad a\ot c \mapsto \sum a\zz{- 1}c\ot a\zz 0, 
\end{equation}
is a left-left weak entwining map.

In particular, if $H$ is a weak bialgebra, then $H$ itself is a right $H$-module 
coalgebra with the action given by the product. Hence, for any right 
$H$-comodule algebra $A$, the map
\begin{equation}\label{Hopf.ent.r}
\psi_R: H\ot A\to A\ot H, \quad h\ot a \mapsto \sum a\zz 0 \ot ha\zz 1
\end{equation}
is a right-right weak entwining map. Therefore, there is a corresponding projection
\begin{equation}\label{pr.Hopf}
p_R: A \otimes H \to A \otimes H, \qquad a\otimes h\mapsto \sum a 1 \zz 0
\otimes h 1 \zz 1,
\end{equation}
and the $A$-coring $\cE= \mathrm{Im} \ p_R $ as in Proposition~\ref{prop.ent.cor}. Explicitly, the structure maps
for $\cE$ come out as $$
\sum b'(a 1 \zz 0 \otimes h 1 \zz 1) b = \sum b' a b \zz 0 \otimes h b \zz 1,
$$ $$ \DE(\sum 1\zz 0 \otimes h 1\zz 1) = \sum (1 \zz 0 \otimes h \sw 1 1 \zz 1)
\otimes_{A} (1 \otimes h \sw 2), $$ $$ \eE (\sum 1 \zz 0 \otimes h 1 \zz
1) = \sum 1 \zz 0 \eps (h 1 \zz 1).
$$
The category of right $\cE$-comodules is isomorphic to the category of
weak relative Hopf modules. In particular, $A$ is such a module, hence
$g= \varrho(1) = \sum 1\zz 0\ot 1\zz 1$ is a grouplike element in $\cE$.
Let $B=A^{co H} = A^{co\cE}_g$. Following \cite{CaeDeG:Gal},  the extension of 
algebras $B\subseteq A$ is called a {\em weak Hopf-Galois $H$-extension} if 
$A$ is a Galois comodule, i.e., $\cE$ is a Galois coring, or, equivalently,
the map 
$$
\mathrm{can}_A:A \otimes _{B} A \to \cE, \qquad \mathrm{can}_A(a
\otimes _{B} b)=a(\sum 1\zz 0 \otimes 1\zz 1) b = \sum ab \zz 0 \otimes b
\zz 1
$$
is an isomorphism of $A$-corings. A weak Hopf algebra $H$ is a weak Hopf-Galois
extension of its (right) coinvariant subalgebra $H^{co H} = \mathrm{Im}\; \pbr$ (cf.\ 
\cite[Proposition~2.7]{CaeDeG:Gal}). 

Obviously, a weak Hopf-Galois extension is
a weak coalgebra-Galois extension. In fact, the above example of a weak Hopf algebra as a weak Hopf-Galois extension is a special case of the general construction of extensions over comodule subalgebras. This construction follows a similar pattern as
in the case of  Hopf algebras (cf.\ \cite[34.2]{BrzWis:cor}) and is described in the 
following 
\begin{example}
Let $H$ be a $k$-flat weak
Hopf algebra, and let $A$ be a left comodule subalgebra of $H$, i.e., a (unital) subalgebra such that $\Delta(A)\subseteq H\ot A$. Let $A^R = \Im\pbr \cap A$. 
Then $J = A^RH$ is a coideal in $H$, since, for all $a\in A^R$ and $h\in H$,
\begin{eqnarray*}
\Delta(ah) &=& \sum a\sw 1h\sw 1\ot a\sw 2h\sw 2 \\
&=& \sum a\sw 1h\sw 1\ot (a\sw 2- \pbr(a\sw 2))h\sw 2  +  \sum a\sw 1h\sw 1\ot 
 \pbr(a\sw 2)h\sw 2 \\
 &=& \sum a\sw 1h\sw 1\ot (a\sw 2- \pbr(a\sw 2))h\sw 2  +  \sum a1\sw 1h\sw 1\ot 
 1\sw 2h\sw 2 \\
 &=&  \sum a\sw 1h\sw 1\ot (a\sw 2- \pbr(a\sw 2))h\sw 2  +  \sum ah\sw 1\ot 
 h\sw 2 \in H\ot J \oplus J\ot H,
 \end{eqnarray*}
 where the third equality follows by \eqref{comal.pr} (remember that $H$ is a right
 $H$-comodule algebra) and the final inclusion is inferred from the assumption $\Delta(a)\in H\ot A$ and from the fact that $\pbr$ is a projection. Furthermore,
 $$
 \eps(ah) = \sum \eps(a1\sw 1)\eps(1\sw 2h) = \eps(\pbr(a)h) = 0.
 $$
 Thus $C = H/J$ is a coalgebra with the coproduct and counit induced by the
 canonical projection $\pi: H\to C$, 
 $\Delta_C = (\pi\ot\pi)\circ\Delta$, $\eps_C \circ \pi =\eps$. Obviously, $C$ is a right $H$-module with the
 multiplication, for all $h\in H$ and $c\in C$,
 $$
 ch = \pi(\th h), \qquad \th\in\pi^{-1}(c).
 $$
In fact, $C$ is a right $H$-module coalgebra, as  $\pi$ is a coalgebra map and $\Delta$ is comultiplicative, hence, for all $h\in H$, $c\in C$ and $\th\in\pi^{-1}(c)$,
\begin{eqnarray*}
\Delta_C(ch) &=& \sum \pi(\th\sw 1h\sw 1)\ot \pi(\th\sw 2h\sw 2) = 
\sum \pi(\th\sw 1)h\sw 1\ot \pi(\th\sw 2)h\sw 2 \\
&=& \sum \pi(\th)\sw 1h\sw 1\ot \pi(\th)\sw 2h\sw 2 =
 \sum c\sw 1h\sw 1\ot c\sw 2h\sw 2.
 \end{eqnarray*}
 This calculation has exactly the same form as a corresponding calculation for
 Hopf algebras. A slightly different computation that makes use of \eqref{eps} and 
 the counitality of $\pi$,  proves the condition \eqref{modco.counit}. Explicitly,
 \begin{eqnarray*}
 \eps_C (chk) &=& \eps_C(\pi(\th hk)) = \eps(\th hk) = \sum\eps(\th h\sw 2)\eps(h\sw 1
 k)\\
 &=& \sum\eps_C(\pi(\th h\sw 2))\eps(h\sw 1
 k) = \sum \eps_C(ch\sw 2)\eps(h\sw 1
 k).
 \end{eqnarray*}
 Thus $(H,C,\re)$,
 where
\begin{equation}\label{re.homog}
 \re :C\ot H\to H\ot C, \qquad c\ot h\mapsto\sum h\sw 1\ot\pi(\th h\sw 2), \qquad 
 \th \in \pi^{-1}(c),
\end{equation}
 is a right-right
 weak entwining structure. Furthermore, $H$ is an entwined module over $(H,C,\re)$ 
 with the
 coaction
 $$
 \varrho^H: H\to H\ot C, \qquad h\mapsto \sum h\sw 1\ot \pi(h\sw 2).
 $$
 Note also that $A\subseteq B = H^{coC}$, as, for all $a\in A$,
 \begin{eqnarray*}
 \varrho^H(a) &=& \sum a\sw 1\ot \pi(a\sw 2) \\
 &= &\sum  a\sw 1\ot \pi(a\sw 2-\pbr(a\sw 2))
 +  \sum  a\sw 1\ot \pi(\pbr(a\sw 2)) = \sum a1\sw 1\ot \pi(1\sw 2),
 \end{eqnarray*}
 where the final equality follows form the facts that $\Delta(a)\in H\ot A$, the definition
 of $\pi$ and by \eqref{comal.pr}. Now define the map
 $$
 \sigma : H\ot C\to H\ot_B H, \qquad h\ot c\mapsto hS\th\sw 1\ot_B \th\sw 2, \qquad 
  \th \in \pi^{-1}(c).
 $$
The map $\sigma$ is well defined, since, if $c=0$, then $\th = \sum_ia^ih^i$, with $a^i\in A^R$, hence
\begin{eqnarray*}
\sum_i h S(a^ih^i)\sw 1\ot_B (a^ih^i)\sw 2 &=& \sum_i h S h^i
\sw 1Sa^i\sw 1\ot_B a^i\sw 2h^i\sw 2\\
&=&\sum_i h S h^i
\sw 1Sa^i\sw 1a^i\sw 2\ot_B h^i\sw 2\\
&=&\sum_i h S h^i
\sw 1\pr(a^i)\ot_B h^i\sw 2\\
&=&\sum_i h S h^i
\sw 1S\pbr(a^i)\ot_B h^i\sw 2=0.
\end{eqnarray*}
The first equality follows by the fact that $S$ is an anti-algebra map, the second one
is a consequence of $\Delta(a)\in H\ot A\subseteq H\ot B$, then the definition of a counit
is used, and finally a relationship between barred and unbarred right projection (cf.\ 
\cite[36.11]{BrzWis:cor} or equation \eqref{s.pbar} below) yields the penultimate
equality. The map $\sigma$ is clearly
a left $H$-module, right $C$-comodule map. Furthermore, for all $h,\th\in H$,
\begin{eqnarray*}
 \sigma \circ \cocan (h\ot_B\th) &=&\sum h\th\sw 1S\th\sw 2\ot_B \th\sw 3 = 
 \sum h\pl(\th\sw 1)\ot_B\th\sw 2\\
 &=&  \sum h\eps(1\sw 1\th\sw 1)1\sw 2\ot_B\th\sw 2 =  \sum h\eps(1\sw 1\th\sw 1)\ot_B1\sw 2\th\sw 2\\
 &=& h\ot\th,
 \end{eqnarray*}
 where the penultimate equality follows by the fact that $A$ is assumed to be a unital
 comodule subalgebra of $H$, so $\Delta(1)\in H\ot A\subseteq H\ot B$. 
 
 In this way we have proven that every comodule subalgebra of a weak Hopf algebra
 $H$ yields a weak coalgebra-Galois extension $B\subseteq H$. Note that
 the canonical entwining map computed directly from the retraction $\sigma$ above (cf.\ Example~\ref{ex.can})
 comes out as, for all $c\in C$, $h\in H$ and $\th\in\pi^{-1}(c)$,
 \begin{eqnarray*}
 \psi^\sigma_R(c\ot h) &=& \cocan(\sum S\th\sw 1\ot_B \th\sw 2h) = \sum S\th\sw 1\th\sw 2h\sw 1\ot \pi(\th\sw 3h\sw 2) \\
 &=& \sum \pr(\th\sw 1)h\sw 1\ot \pi(\th\sw 2h\sw 2) =  \sum  1\sw 1h\sw 1\ot \pi(\th1\sw 2h\sw 2) \\
 &=&  \sum  h\sw 1\ot \pi(\th h\sw 2),
 \end{eqnarray*}
where  the penultimate equality follows by the left comodule version of \eqref{comal.pl}
($H$ is a left $H$-comodule algebra). Thus the canonical entwining $\psi^\sigma_R$ coincides with the entwining $\re$ in \eqref{re.homog}, as expected.
\end{example}

\section{A quest for an invertible weak entwining structure}
Recent generalisations of the Schneider Theorem I (cf.\ \cite{Sch:pri}) 
and the Kreimer-Takeuchi theorem (cf. \cite{KreTak:hop}) to 
general  coalgebra-Galois extensions and coring-Galois extensions use 
bijectivity of the canonical entwining structure (cf.\ \cite{Brz:gal}, \cite{SchSch:gen}) or 
an entwining structure
over a non-commutative ring (cf.\ \cite{Boh:Gal}).  As the aim of the present paper is to extend these results
to weak coalgebra-Galois extensions, we need to find a proper definition of an
invertible weak entwining structure. If we make the obvious choice, however, then
our quest suffers 
immediate setback 
because of the following
\begin{lemma}\label{lemma.bij}
If $(A,C, \psi_R)$ is a right-right weak entwining structure such that $\psi_R$
is a bijective map, then $(A,C, \psi_R)$ is an entwining structure.
\end{lemma}
\begin{proof}
For the $k$-linear inverse of $\psi_R$, write 
$\psi_R^{-1}=\sumA c_E \otimes a^E$, so that, for all $a\in A$ and $c\in C$,
$$
c\otimes a = \sum_{\alpha, E} {c\spa}_{E} \ot {a\sba}^{E}, \qquad a\otimes c = \sum_{\alpha, E} {a^{E}}_{\alpha} \ot {c_{E}}^{\alpha} .$$
Applying $\psi_R^{-1}$ to relation (\ref{re1})  we obtain, 
for all $a,  \tilde{a}\in A$ and $c\in C$,
$$ \sumAB c_{EF} \otimes a^F
\tilde{a}^E=\sumA c_E \otimes (a \tilde{a})^E.
$$
With these at hand, take any $a\in A$ and $c\in C$, and compute 
\begin{eqnarray*}
 \suma a1\sba \otimes c\spa
     &=& \suma \psi_R \circ \psi_R^{-1}(a1\sba \otimes c\spa)\\
 & =& \sum_{\alpha, E,F} \psi_R ({c\spa}_{EF} \otimes
a^F {1\sba}^E)
     =\sum_F \psi_R (c_F \otimes a^F)= a \otimes c .
 \end{eqnarray*}
Now taking $a=1$, we deduce that 
 $\suma 1\sba \otimes c\spa = 1 \otimes
c$.  This equality, together with the  relation (\ref{re2}) imply
 $$\suma a\sba \eps_C (c \spa)=\suma 1\sba \eps_C(c \spa) a = \eps_C (c) a,
 $$
  i.e.,
 $(A,C,\psi_R)$ is a right-right entwining structure as claimed.
\end{proof}
 
 Lemma~\ref{lemma.bij} means that, if one wants to deal with weak entwining
 structures, one cannot assume that the entwining map be bijective. Instead we 
 propose
\begin{definition}\label{def.inv}
An {\em invertible weak entwining structure} is a quadruple $(A,C, \re, \le)$ such that 
 \begin{blist}
 \item $(A,C,\re)$ is a right-right weak entwining structure and $(A,C,\le)$ is a 
 left-left weak entwining structure;
 \item $
\re \circ \le = p_R$ and $ \le \circ \re = p_L$; 
 \item for all $c\in C$, 
 $$\sumA\eps_C(c_E)1^E=\suma 1\sba \eps_C (c\spa
).
 $$
 \end{blist}
\end{definition}
The introduction of a bijective entwining map in all generalisations of Schneider's and 
the Kreimer-Takeuchi theorems, replaces the original assumption that a Hopf algebra
has a bijective antipode. Indeed,  the entwining map associated to a Doi-Koppinen 
datum $(A,C,H)$ is bijective provided $H$ has a bijective antipode. The notion of 
an invertible weak entwining structure introduced in  Definition~\ref{def.inv} is 
motivated by the following observation.
\begin{proposition}\label{prop.doi.inv}
Let $H$ be a weak Hopf algebra with a bijective antipode, 
$A$ a right $H$-comodule algebra and $C$ a right 
$H$-module coalgebra. Then $(A,C,\re,\le)$  with $\re$ given by equation 
(\ref{doi.ent.r}) and 
$$
\le : A\ot C\to C\ot A, \qquad a \otimes c\mapsto  \sum c S^{-1} a\zz 1 \otimes a \zz 0,
$$
is an invertible weak entwining structure.
\end{proposition}
\begin{proof}
First we need to prove that $(A,C,\le)$ is a left-left weak entwining structure. 
Suffices it to check that  there exists a weak bialgebra $\bar{H}$ such that $A$
is a left $\bar{H}$-comodule algebra, $C$ is a left $\bar{H}$-module coalgebra
and  the map $\le$ has the form given in equation (\ref{doi.ent.l}). Take $\bar{H}=H^{op}$, 
the opposite algebra to $H$. When viewed as an element of $\bar{H}$, an element $h$ of $H$ is denoted by $\bar{h}$. The right $H$-multiplication on $C$, makes $C$ a left $\bar{H}$-module by the formula
$$
{\bar h}c  = ch.
$$
Since the right action of $H$ on $C$ is comultiplicative, so is the derived left 
$\bar{H}$-action.
Furthermore, condition (\ref{modco.counit}) implies for all $h,k\in H$ and $c\in C$,
$$
\eps_C (\bar{k}\bar{h}c) = \eps_C (chk) = \sum \eps_C(ch\sw 2)\eps(h\sw 1k) =\sum \eps(\bar{k}\bar{h}\sw 1)\eps_C(\bar{h}\sw 2c) ,
$$
(note the use of equation (\ref{modco.counit})). This is one of the equivalent conditions for $C$ to be a left $\bar{H}$-module coalgebra.
Consider a $k$-linear map
$$
\lambda : A\to \bar{H}\ot A, \qquad a\mapsto \sum S^{-1}a\zz 1\ot a\zz 0.
$$
Since the antipode is an anti-coalgebra map, so is $S^{-1}$, hence $\lambda$ is
a left $\bar{H}$-coaction. Furthermore, $\varrho$ is a multiplicative map, while $S^{-1}$
is an anti-algebra map, thus  $\lambda : A\to \bar{H}\ot A$ is a multiplicative map. 
Finally, we can use equations (\ref{s.pi}) and (\ref{comal.pl.one}) to compute
$$
\sum \pr (S^{-1}1\zz 1)\ot 1\zz 0 = \sum S^{-1}\pl (1\zz 1)\ot 1\zz 0 = \sum S^{-1} 1\zz 1\ot 1\zz 0.
$$
This is the left-handed version of condition  (\ref{comal.pl.one}). Thus we conclude 
that $A$ is a left $\bar{H}$-comodule algebra, $C$ is a left $\bar{H}$-module coalgebra
and $\le$ has the form (\ref{doi.ent.l}) with the coaction $\lambda$. Hence $(A,C,\le)$ is 
a left-left weak entwining structure.

Next we prove that the conditions (b) in Definition~\ref{def.inv}
are satisfied. The projections corresponding to weak entwining 
maps $\re$ and $\le$ come
out as
$$p_R(a
\otimes c)= \sum a1 \zz 0 \otimes c 1 \zz 1, \qquad 
p_L(c \otimes a)=\sum c S^{-1} 1 \zz 1 \otimes 1 \zz 0 a.
$$
Take any $a\in A$ and $c\in C$ and compute
\begin{eqnarray*}
\le \circ \re (c \otimes a)
    &=& \le (\sum a \zz 0 \otimes c a \zz 1)=\sum ca \zz 1\sw  2 S^{-1} a \zz 1\sw 1 \otimes a
\zz 0 \\
    &=& \sum c S^{-1}( a \zz 1\sw 1 S a \zz 1\sw 2) \otimes a \zz 0  
    = \sum c S^{-1} \pl (a \zz 1) \otimes a \zz 0  \\
    &=& \sum cS^{-1} 1\zz 1 \otimes 1\zz 0 a 
    = p_L(c \otimes a).
\end{eqnarray*}
The third equality follows from the fact that the antipode (and hence also its inverse) is
an anti-algebra map, the fourth one is the defining property of the antipode in a weak
Hopf algebra. The fifth equality follows from the definition of a right 
$H$-comodule algebra, equation (\ref{comal.pl}).  
Next use the first of equations (\ref{s.pi}) and the first of equations (\ref{pi.s}) to
 note that, for all $h\in H$,
\begin{gather}\label{s.pbar}
S^{-1}\pr(h)=\pl(S^{-1}h)=\sum \eps(h1\sw 1 )1 \sw 2 =\pbr(h).
\end{gather}
Use this to compute
\begin{eqnarray*}
\re \circ \le (a \otimes c)
    &=& \re (\sum cS^{-1}a\zz 1 \otimes a \zz 0)=\sum a \zz 0 \otimes cS^{-1}(a \zz 1\sw  2)
a \zz 1 \sw 1 \\
    &=& \sum a \zz 0 \otimes c S^{-1}((Sa\zz 1\sw  1)a \zz 1\sw 2)
    = \sum a \zz 0 \otimes c S^{-1} \pr(a \zz 1)\\
    &=& \sum a \zz 0 \otimes c \pbr (a \zz 1) 
    = \sum a1\zz 0 \otimes c 1 \zz 1 
    = p_R(a \otimes c).
\end{eqnarray*}
 Again the third equality follows from the fact that the antipode (and hence also its inverse) is
an anti-algebra map and  the fourth one is the defining property of the antipode in a weak Hopf algebra. The fifth equality follows from property (\ref{s.pbar}), while the sixth one is 
the consequence of the fact that $A$ is a right $H$-comodule algebra, hence 
equation (\ref{comal.pr}) holds.

Finally we need to check that the maps $\le$ and $\re$ satisfy property (c) in
Definition~\ref{def.inv}, i.e., that, for all $c\in C$,
$\sum 1 \zz 0 \eps_C (c1\zz 1) = \sum \eps_C (c S^{-1}1\zz 1
)1\zz 0$. Compute
\begin{eqnarray*}
\sum \eps_C(c S^{-1}1\zz 1) 1\zz0
    &=& \sum \eps_C(c1\sw 2)\eps(1\sw 1 S^{-1}1\zz 1)1\zz 0\\
&=&  \sum \eps_C(c1\sw 2\eps(1\sw 1 S^{-1}1\zz 1))1\zz 0\\
&=&  \sum \eps_C(c\pl( S^{-1}1\zz 1))1\zz 0\\
&=&  \sum \eps_C(c\pbr(1\zz 1))1\zz 0 = \sum \eps_C (c 1\zz 1)1\zz 0.
\end{eqnarray*}
The first equality follows from the compatibility between the counit and the action of
$H$ on a right module coalgebra $C$, equation~(\ref{modco.counit}). The third equality
is simply the definition of $\pl$, while the fourth equality follows from property (\ref{s.pbar}). 
Finally, the last equality is obtained by the use of equation~(\ref{comal.pr.one}). 
Thus we conclude that $(A,C,\re,\le)$ is an invertible weak entwining structure
as claimed.
\end{proof}

In particular we obtain
\begin{corollary}\label{cor.doi.inv}
Let $H$ be a weak Hopf algebra with a bijective antipode and let
$A$ be a right $H$-comodule algebra. Then $(A,H,\re,\le)$  with 
$\re$ given by equation 
(\ref{Hopf.ent.r}) and 
\begin{equation}\label{Hopf.ent.l}
\le : A\ot C\to C\ot A, \qquad a \otimes h\mapsto  \sum h S^{-1} a\zz 1 \otimes a \zz 0,
\end{equation}
is an invertible weak entwining structure.
\end{corollary}
\begin{proof}
This follows immediately from the fact that $H$ itself is a right $H$-module 
coalgebra with the action given by the product, and from Proposition~\ref{prop.doi.inv}.
\end{proof}

As an immediate consequence of Definition~\ref{def.inv} one obtains the following
\begin{lemma}\label{lemma.inv}
Let  $(A,C, \re, \le)$ be an invertible weak entwining structure. Then:
\begin{zlist}
\item $\re\circ p_L = \re$,
\item $\le\circ p_R = \le$.
\end{zlist}
\end{lemma}
\begin{proof} (1) Note that $\re(C\ot A)\subseteq \mathrm{Im}\ p_R$. Since $p_R$ is a
projection, $p_R\circ \re = \re$, and the condition (a) in Definition~\ref{def.inv} implies
$$
\re\circ p_L = \re\circ\le\circ\re = p_R\circ\re = \re.
$$
Statement (2) is proven in a similar way.
\end{proof}

As is often the case with (weak) entwining structures, the notion of an invertible
weak entwining structure has the clearest meaning in terms of corings.
\begin{proposition}\label{prop.inv.cor}
Let $(A,C,\re,\le)$ be an invertible weak entwining structure and let
$\cC=\mathrm{Im} \ p_R$ and $\cD = \mathrm{Im} \ p_L$ be the corresponding
$A$-corings.  Then the restrictions of the entwining maps
$$
\le: \cC \to \cD, \qquad \re : \cD \to \cC
$$
are inverse isomorphisms of $A$-corings.
\end{proposition}
\begin{proof}
Since $p_R$ and $p_L$ are projections, the conditions (b) in Definition~\ref{def.inv}
imply that restrictions of $\re$ and $\le$ to $\mathrm{Im} \ p_L$ and 
$\mathrm{Im} \ p_R$ respectively,  are inverse 
 isomorphisms of $k$-modules. Furthermore, for all $a\in A$ and $c\in C$,
$$
\re (\sumA a (c_E \otimes 1^E))=\re (\sumAB c_{EF} \otimes a^E1^F)= 
\re (\sumA c_E \otimes a^E) = p_R(a\otimes c)= \suma a 1\sba \otimes c \spa,
$$
where the second equality follows from \eqref{le1} and the 
penultimate equality follows from the property (b) in Definition~\ref{def.inv}.
On the other hand
$$
a \re (\sumA c_E \otimes 1^E) = a  \re (\le (1 \otimes c))= a
 p_R (1\otimes c) = \suma a1\sba \otimes c \spa,
$$
where  again the second equality 
follows from the property (b) in Definition~\ref{def.inv}. Thus $\re$ is a left $A$-module
map. 
Moreover, Lemma~\ref{lemma.inv} implies
$$
\re(\sumA c_E \otimes 1^E a) =\re\circ p_L (c\ot a)= \re(c\otimes a) = \suma a\sba \otimes c \spa .
$$
On the other hand
$$
\re (\sumA c_E \otimes 1^E)  a = \suma (1 \sba \otimes c \spa)  a = \sumab 1\sba a\sbb \otimes c^{\alpha\beta} = 
\suma a \sba \otimes c \spa,
$$
where the last equality follows from equation (\ref{re1}). 
Hence $\re$ is an $(A,A)$-bimodule map.  Similarly it can be shown that
$\le$ is an $(A,A)$-bimodule map. Thus $\re$ is an $(A,A)$-bimodule isomorphism
with the inverse $\le$.

Next we show that the map  $\re$ is counital. Take any $c\in C$ and $a\in A$ 
and compute
\begin{eqnarray*}
\eC(\re ( \sumA c_E \otimes 1^Ea)) &=& \suma a\sba \eps_C (c\spa) 
= \suma 1\sba \eps_C(c
\spa)a\\
& =& \sumA \eps_C (c_E)1^E a = \eD (\sumA c_E\otimes 1^Ea).
\end{eqnarray*}
To derive the first equality we have used Lemma~\ref{lemma.inv}(1) and then the definition of the counit in $\cC$. The second
equality follows from equation (\ref{re2}), while the third one is the consequence of
propety (c) in Definition~\ref{def.inv}. The final equality follows from the definition of
the counit in $\cD$. Thus $\re$ is a counital map. 
Similarly one shows that $\le$ is a counital map. 

Finally we need to prove that $\re$  and $\le$ are comultiplicative maps. Take
any $a\in A$ and $c\in C$, and use Lemma~\ref{lemma.inv}(1) and the definition of the coproduct in $\cC$ to compute
$$
\DC \circ \re(\sumA c_E \otimes 1^E a)=\DC(\suma a\sba \otimes c \spa)= \suma a\sba
\otimes {c\spa} \sw 1 \otimes {c\spa} \sw 2. $$
On the other hand,
\begin{eqnarray*}
(\re \otimes _A \re)\circ \DD(\sumA c_E \otimes 1^Ea)
    &=& \sumA \re(c_E \sw 1 \otimes 1) \otimes _A \re (c_E\sw 2 \otimes
    1^Ea)\\
    &=& \sum_{E,F} \re (c \sw 1 _E \otimes 1) \otimes _A \re ( c \sw 2 _F
    \otimes 1^{EF} a)\\
 &=& \sum_{E} \re (c \sw 1 _E \otimes 1^E) \otimes _A \re ( c \sw 2 
    \otimes a)\\
    &=& \sumab 1  \sba \otimes c \sw 1 \spa \otimes _A a \sbb \otimes c
    \sw 2 \spb\\
    &=& \sum_{\alpha,\beta,\gamma}1\sba a _{\beta \gamma} \otimes c \sw 1 ^{\alpha \gamma} \ot_A 1\otimes c
    \sw 2 \spb\\
    &=& \suma a \sba \otimes c \spa \sw 1  \ot_A 1\otimes c  \spa \sw 2 = \DC (\suma  a
    \sba \otimes c \spa ).
\end{eqnarray*}
The first equality is the definition of the coproduct in $\cD$, the second follows from the
 definition of a left-left entwining structure, more precisely, equation~(\ref{le3}). Next we
 use the definition of the left $A$-multiplication in $\cD$ and the fact that $\re$ is an
 $(A,A)$-bimodule map. The property (b) in Definition~\ref{def.inv}
   implies the fourth equality.  Then we use
 the definition of right $A$-multiplication in $\cC$. The final two equalities
 follow from the properties of a right-right entwining structure, equations (\ref{re1}) and (\ref{re3}). Similarly one proves that $\le$ is comultiplicative and thus completes the proof
 of the proposition.
\end{proof}
 
In the case of entwining structures, the inverse $\psi^{-1}$ 
of a right-right entwining map $\psi$ is a 
left-left entwining map. Furthermore,  if $A$ is a right entwined module over $(A,C,\psi)$,
 then it is a left entwined module over $(A,C,\psi^{-1})$ (cf.\ \cite[Section~6]{Brz:mod}).
 Similarly, for an invertible weak entwining structure one proves 
\begin{corollary}\label{cor.left.coa}
Let $(A, C, \re,\le)$ be an invertible weak entwining structure.
If $A \in \rem$, then $A \in \lem$ with the coaction
$$ ^{A}\!\varrho (a) = \le (\sum a1\zz 0 \otimes 1\zz 1).$$
\end{corollary}
\begin{proof}
If $A \in \rem$, 
 then $A$ is a  right $\cC$-comodule.  The corresponding grouplike element is $g= \varrho^A(1) =\sum 1 \zz 0 \otimes 1 \zz 1$. But then $A$ is  also a left $\cC$-comodule with the coaction $a\mapsto ag\ot_A 1$. By Proposition~\ref{prop.inv.cor}, the map
  $\le :\cC\to \cD$ is an isomorphism of
 $A$-corings, hence the left $\cC$-coaction of $A$ induces a left $\cD$-coaction on $A$,
 $$
{^A \varrho}(a)=\le (ag) \otimes _A 1 \simeq \le (ag).
$$
Thus $A \in \lem$ with coaction ${^A \varrho (a)} = \le (\sum a 1\zz 0
\otimes 1 \zz 1)$, as stated.
\end{proof}

\section{Weak coalgebra Galois extensions with coseparable coalgebras}
Once we have understood how to define an invertible weak entwining structure
we can proceed to apply Theorem~\ref{thm.mod.main} to deduce a criterion for 
an algebra to be a weak coalgebra-Galois extension. This is a subject of the present
and following sections.

First recall that a coalgebra $C$ is called a {\em coseparable}  coalgebra provided  
the coproduct has a retraction in the category of $C$-bicomodules.
Equivalently, $C$ is a coseparable coalgebra if there exists a {\em
cointegral}, i.e., a $k$-module map $\delta :C\ot C\to k$ 
that is {\em colinear}, meaning,
 for all $c,c'\in C$,
\begin{equation}\label{colin}
\sum c\sw 1\delta(c\sw 2\ot c') = \sum \delta(c\ot c'\sw 1)c'\sw 2,
\end{equation}
and such
that   $\delta\circ\Delta_C = \eps_C$.
  Over an algebraically closed field the notion of   coseparability is equivalent to
the notion of cosimplicity. 
\begin{theorem} \label{thm.main} 
Let  $(A,C,\re,\le)$ be an invertible weak
entwining structure with a $k$-projective algebra $A$ and $k$-projective
 coalgebra $C$. Denote by $p_R$ the projection corresponding  to 
 $(A,C,\re)$ and 
 by $\cC =\mathrm{Im}\ p_R$ the corresponding $A$-coring. Suppose that
\begin{zlist} 
\item $A$ is a right weak entwined module with product $\mu$ and coaction
$\varrho^A$; 
\item the map ${\widetilde\can}_A:A \ot A \to\cC$, $a\ot a'\to a\varrho^A(a')$ is 
 surjective;  
\item $C$ is a  coseparable coalgebra.
\end{zlist} 
Then $A$ is a weak $C$-Galois extension of the coinvariants $B$ and is
$C$-equivariantly projective as a left $B$-module (i.e., $A$ is a projective
 left $B$ module and the multiplication map
 $B\ot A\to A$ has a left $B$-linear, right $C$-colinear section). 
\end{theorem}

The strategy of the proof of Theorem~\ref{thm.main} is to show that it is a special case
 of Theorem~\ref{thm.mod.main}. As the first step we need to show that 
 there is a left $\cC$-comodule map $\cC\to A\ot A$ that splits the canonical map $\widetilde{\can}_A: A\ot A\to \cC$. 
By Corollary~\ref{cor.left.coa}, given an invertible weak entwining structure 
such that $A$ is a right
 weak entwined module, $A$ is a left weak entwined module.  In the case of 
 invertible entwining structures, there is a bijective correspondence 
 between left $\cC$-colinear maps $\cC = A\otimes C\to A\ot A$ with
 left $C$-colinear maps $C\to A\ot A$. A similar statement can be proven for
 invertible weak entwining structures.
 
 \begin{lemma} \label{lemma3} 
Given an invertible weak entwining structure $(A,C,\re,\le)$ with 
$A \in \rem$, view $A\ot A$ as   a left $\cC = (\mathrm{Im}\ p_{R})$-comodule with
 the coaction 
 $${}^{A\ot A}\!\varrho :A\ot A\to \cC\ot A, \qquad 
 a\ot a'\mapsto \sum a1\zz 0\ot 1\zz1\ot a'.$$
  Then there is the bijective 
 correspondence between left
$\cC$-colinear maps $f:\cC \to A \ot A$ and $k$-linear maps
$\hat{f}: C \to A \ot A$ such that
\begin{equation} \label{eq.colin}
(^{A}\!\varrho \ot A) \circ \hat{f} = (p_{L} \ot A) \circ (C \ot 
\hat{f}) \circ \Delta_C,
\end{equation}
 where ${}^{A}\!\varrho : A\to C\ot A $ is a left $C$-coaction as in 
 Corollary~\ref{cor.left.coa}.
\end{lemma}
\begin{proof}
For all $c\in C$, write $\hat{f}(c) = \sum c\sv{1} \ot 
c\sv{2}$ and apply $\re \ot A$ to equation~\eqref{eq.colin} to obtain
\begin{equation}\label{first.step}
\sum\re \circ \le( c\sv{1} 1\zz{0} \ot 1\zz{1}) \ot c\sv{2} = \sum \re \circ 
p_{L}(c\sw{1} \ot c\sw{2}\sv{1}) \ot c\sw{2}\sv{2}.
\end{equation}
Since, for all $a\in A$,  $\sum a1\zz{0} \ot 1\zz{1} \in \mathrm{Im}\ p_{R}$, 
 the definition of an invertible
 weak entwining structure, Definition~\ref{def.inv}(b), implies that  $\re \circ 
\le(\sum a1\zz{0} \ot 1\zz{1}) = \sum a1\zz{0} \ot 1\zz{1}$. 
Applying  Lemma~\ref{lemma.inv} to the right hand side of \eqref{first.step},
one obtains that the equation  \eqref{eq.colin} is equivalent to
\begin{equation}\label{second.step}
\sum c\sv{1} 1\zz{0} \ot 1\zz{1} \ot c\sv{2} = \sum \re(c\sw{1} \ot 
c\sw{2}\sv{1}) \ot c\sw{2}\sv{2},
\end{equation}
i.e.,
\begin{equation}\label{third.step}
\sum c\sv{1} 1\zz{0} \ot 1\zz{1} \ot_{A} 1 \ot c\sv{2}
= \sum (1 \ot c\sw{1}) \ot_{A} c\sw{2}\sv{1} \ot c\sw{2}\sv{2}.
\end{equation}

Take a $k$-linear map $\hat{f}: C\to A\ot A$ that satisfies condition \eqref{eq.colin} and define a left $A$-module map $f: \cC \to A\ot A$ by $f(\suma a1\sba \ot c\spa) = \suma a1\sba \hat{f}(c\spa)$. Note that,
\begin{eqnarray*}
^{A \ot A}\!\varrho \circ f(\suma a1\sba \ot c\spa)
&=& \suma a1\sba c\spa\sv{1} 1\zz{0} \ot 1\zz{1} \ot_{A} 1 \ot 
c\spa\sv{2} \\
&=& \suma a1\sba (1 \ot c\spa\sw{1}) \ot_{A} c\spa\sw{2}\sv{1} \ot 
c\spa\sw{2}\sv{2} \\
&=& \suma a1\sba \ot c\spa\sw{1} \ot_{A} c\spa\sw{2}\sv{1} \ot 
c\spa\sw{2}\sv{2},
\end{eqnarray*}
where the second equality follows from property \eqref{second.step}. 
On the other hand
\begin{eqnarray*}
(\cC\ot_A f)\circ\DC (\suma a1\sba \ot c\spa) 
&=&  \sumab a1\sba \ot c\sw{1}\spa \ot_{A} f(1_\beta \ot c\sw{2}^\beta) \\
 &=& \sum_{\alpha,\beta} a1\sba
  \ot c^{\alpha}\sw{1} \ot_{A} 1_\beta c^{\alpha}\sw{2}^\beta\sv{1} \ot 
c^{\alpha}\sw{2}^\beta\sv{2} \\
&=& \sum_{\alpha,\beta,\gamma} a1\sba1_{\beta\gamma}
  \ot c^{\alpha}\sw{1}^\gamma \ot_{A} c^{\alpha}\sw{2}^\beta\sv{1} \ot 
c^{\alpha}\sw{2}^\beta\sv{2} \\
&=& \sumab a1\sba1_\beta  \ot c^{\alpha\beta}\sw{1}
 \ot_{A} c^{\alpha\beta}\sw{2}\sv{1} \ot 
c^{\alpha\beta}\sw{2}\sv{2} \\
 &=& \suma a1\sba \ot c\spa\sw{1} \ot_{A} c\spa\sw{2}\sv{1} \ot 
c\spa\sw{2}\sv{2}.
\end{eqnarray*}
Here we first use the definition of the coproduct in $\cC$ in Proposition~\ref{prop.ent.cor}, then the definition of $f$ in terms 
of the map $\hat{f}$, next the definition of right $A$-multiplication in $\cC$, and 
finally properties \eqref{re1} and \eqref{re3} of a right-right
weak entwining structure. This proves that $f$ is a left $\cC$-colinear map.

Conversely, take a left $\cC$-colinear map $f:\mathrm{Im}\ p_{R} \to A 
\ot A$ and define $\hat{f}: C\to A\ot A$ by $\hat{f}(c) = f(\suma 1\sba \ot c\spa)$. Then
\begin{eqnarray*}
    \sum c\sv{1} 1\zz{0} \ot 1\zz{1} \ot_{A} 1 \ot c\sv{2}
    &=& \suma  {}^{A \ot A}\!\varrho(f(1\sba \ot c\spa)) \\
    &=& \suma 1\sba \ot c\spa\sw{1} \ot_{A} f(1 \ot c\spa\sw{2}) \\
    &=& \sumab 1\sbab \ot c\sw{1}\spb \ot_{A} f(1 \ot c\sw{2}\spa) \\
    &=& \suma (1 \ot c\sw{1}) \cdot 1\sba \ot_{A} f(1 \ot c\sw{2}\spa) \\
    &=& \sum (1 \ot c\sw{1}) \ot_{A} c\sw{2}\sv{1} \ot c\sw{2}\sv{2},
    \end{eqnarray*}
where the first equality follows from the definition of the left $\cC$-coaction in
$A$ in terms of the grouplike element $\sum 1\zz 0\ot 1\zz 1$. The second equality
is a consequence of the fact that $f$ is a left $\cC$-colinear map, while the third one
is the property (\ref{re3}) of a right-right entwining structure. The fourth equality
follows from the definition of the right $A$-action in $\cC$, and the last one is obtained  
by the use of the fact that $f$ is a left $A$-linear map and then by the definition of
$\hat{f}$.  So we obtain \eqref{third.step}, i.e., the condition \eqref{eq.colin} is 
satisfied, as required.
\end{proof}

Before the start of the proof of Theorem~\ref{thm.main} it is also useful to prove 
the following

\begin{lemma} \label{lemma4}
Take an  invertible weak entwining structure $(A,C,\re,\le)$ with a $k$-flat coalgebra
$C$. Assume that $A$ is a weak entwined module as in assumption (1) in Theorem~\ref{thm.main} and that the canonical map $\widetilde{\can}_A
:A \ot A \to \mathrm{Im}\ p_R$ has a $k$-linear section  $\tau$. Set
$\hat\tau(c)=\tau(\suma 1\sba \otimes c\spa)$. Then 
\begin{equation} \label{eq.coactns}
(C \ot \mu \ot C) \circ (^A\!\varrho \ot \varrho^A) \circ \hat\tau 
= (C \ot \psi_R) \circ (C \ot C \ot 1) \circ \Delta_C,
\end{equation}
 where ${}^{A}\!\varrho : A\to C\ot A $ is a left $C$-coaction as in 
 Corollary~\ref{cor.left.coa}.
\end{lemma}
\begin{proof}
This can be shown as follows. Write $\hat\tau(c) = \sum c\sv{1} \ot 
c\sv{2}$. Since $\tau$ is a section of $\widetilde{\can}_A$,
\begin{equation}\label{eq.section}
\widetilde{\can}_A \circ \tau(\suma 1\sba \ot c\spa):= \sum c\sv{1} c\sv{2}\zz{0} 
\ot c\sv{2}\zz{1} = \suma 1\sba \ot c\spa.
\end{equation}
Start with the identity
\begin{eqnarray*}
\sum c\sv{1}\zz{-1} \ot c\sv{1}\zz{0} c\sv{2}\zz{0} \ot c\sv{2}\zz{1}
&=& \sum \le(c\sv{1} 1\zz{0} \ot 1\zz{1})1\zz{0'} \re(1\zz{1'} \ot 
c\sv{2}) \\
 &=& \sum_{\alpha, E}
1\zz{1}_{E} \ot (c\sv{1} 1\zz{0})^{E} 1\zz{0'} c\sv{2}\sba \ot 
1\zz{1'}\spa,
\end{eqnarray*}
which is simply the definition of left and right $C$-coactions and also
uses \eqref{re1}.
Applying $\re \ot C$ we obtain
\begin{eqnarray*}
    \sum \re(c\sv{1}{}\zz{-1} \ot c\sv{1}{}\zz{0} c\sv{2}{}\zz{0}) \ot 
    c\sv{2}{}\zz{1} \!\!\!
    &=&\!\!\!\!\!\!\! \sum_{\alpha,\beta,\gamma ,\delta, E}
 (c\sv{1} 1\zz{0})^{E}{}_{\beta} 1\zz{0'}{}_{\gamma}
    c\sv{2}{}\sba{}_{\delta} \ot 1\zz{1}{}_{E}{}^{\beta\gamma\delta}
    \ot 1\zz{1'}{}\spa \\
    &=& \!\!\!\! \sum_{\alpha,\beta,\gamma ,\delta}
 c\sv{1} 1\zz{0} 1_{\beta} 1\zz{0'}{}_{\gamma} 
    c\sv{2}{}_{\alpha\delta} \ot 1\zz{1}{}^{\beta\gamma\delta} \ot 
    1\zz{1'}{}\spa \\
    &=& \sum_{\alpha,\delta}c\sv{1} 1\zz{0} (1\zz{0'} c\sv{2}{}_{\alpha})_{\delta} 
    \ot 1\zz{1}{}^{\delta} \ot 1\zz{1'}{}\spa \\
    &=& \sum_{\delta}c\sv{1} 1\zz{0} c\sv{2}{}\zz{0}{}_{\delta} \ot 
    1\zz{1}{}^{\delta} \ot c\sv{2}{}\zz{1} \\
 &=& \sum c\sv{1}  c\sv{2}{}\zz{0} \ot 
    c\sv{2}{}\zz{1}\sw 1 \ot c\sv{2}{}\zz{1}\sw 2 \\
    &=& \suma 1\sba \ot c\spa{}\sw{1} \ot c\spa{}\sw{2} \\
    &=& \sumab 1\sbab \ot c\sw{1}{}\spb \ot c\sw{2}{}\spa = 
 \suma \re (c\sw 1\ot 1_\alpha)\ot c\sw 2^\alpha .
    \end{eqnarray*}
   The first and the third equality follow from the (multiple)
  use of property \eqref{re1} of a right-right weak
 entwining structure and the second follows from Definition~\ref{def.inv}~(b). The 
 fourth and fifth equalities are consequences of the fact that $A$ is a right
 entwined module. The sixth equality follows from equation~\eqref{eq.section}, 
  while the penultimate equality follows from equation \eqref{re3}.
 Now we apply $\le \ot C$ and use the fact that 
    $$\sum c\sv{1}\zz{-1} \ot c\sv{1}\zz{0} c\sv{2}\zz{0} \ot c\sv{2}\zz{1}
    \in \mathrm{Im}\ p_{L}\otimes C,$$
  to compute
    \begin{eqnarray*}
	\sum c\sv{1}\zz{-1} \ot c\sv{1}\zz{0} c\sv{2}\zz{0} \ot c\sv{2}\zz{1}
	&=& \suma \le \circ \re(c\sw{1} \ot 1\sba) \ot c\sw{2}{}\spa \\
	&=& \sum_{\alpha,E}c\sw{1}{}_{E} \ot 1^{E} 1\sba \ot c\sw{2}{}\spa \\
	&=& \sum_{\alpha, E}c\sw{1} \ot \eps_C(c\sw{2}{}_{E}) 1^{E} 1\sba \ot 
	c\sw{3}{}\spa \\
	&=& \sumab c\sw{1} \ot 1\sbb 1\sba \ot c\sw{2}{}\spba \\
	&=& \suma c\sw{1} \ot 1\sba \ot c\sw{2}{}\spa.
	\end{eqnarray*}
	  The second equality follows from the definition of an invertible
  weak entwining structure, Definition~\ref{def.inv}~(b). Then we have used
  property \eqref{le4} of a left-left weak entwining structure, then Definition~\ref{def.inv}~(c) and \eqref{re4}. The last equality follows
  by \eqref{re1}. Thus we obtain \eqref{eq.coactns} as required.
    \end{proof}

We are now ready to prove Theorem~\ref{thm.main}.

\begin{proof}{(Theorem \ref{thm.main})}
Take $\cC=\mathrm{Im}\ p_R$.  By the $k$-projectivity assumptions, $A\ot C$ is a 
projective $k$-module. Since $\cC$ is a direct summand in $A\ot C$, it is a projective
$k$-module too. 
Let $\tau$ be a $k$-linear splitting of $\widetilde{\can}_A: A \ot A \to 
\cC$, and set  $\hat\tau(c) = \tau(\suma 1\sba \ot c\spa)$. 
Let $\delta$ be a cointegral of $C$ and consider 
\begin{equation} 
\hat\kappa = (\delta \ot A \ot A) \circ (C \ot ^A\!\varrho \ot A)
\circ (C \ot \hat\tau) \circ \Delta_C, 
\end{equation} 
where we view $A$ as a left $C$-comodule with coaction $^A\!\varrho$ as in
Corollary~\ref{cor.left.coa}. Using the colinearity of $\delta$, one easily
checks that $\hat{\kappa}$ is a left $C$-colinear map, i.e.,
$$
({}^A\!\varrho\ot A)\circ \hat{\kappa} = (C\ot\hat{\kappa})\circ \Delta_C.
$$
Since ${}^A\!\varrho(A)\subseteq \Im p_L$, and $p_L$ is a 
projection, this implies that  $\hat\kappa$ satisfies 
the property \eqref{eq.colin} in Lemma~\ref{lemma3}. Thus the map
\begin{equation}
\kappa(\suma a 1\sba \ot c\spa) = \suma a 1\sba \hat\kappa(c\spa)
\end{equation}
is a left $\cC$-comodule map. Explicitly, writing $\hat{\tau}(c) = \sum c\sv 1\ot c\sv 2$,
the map $\kappa$ comes out as
\begin{equation}
\kappa(\suma a 1\sba \ot c\spa) = \suma a 1\sba \delta(c\spa{}\sw{1} \ot 
c\spa{}\sw{2}{}\sv{1}{}\zz{-1}) c\spa{}\sw{2}{}\sv{1}{}\zz{0} \ot 
c\spa{}\sw{2}{}\sv{2}.
\end{equation}
Note that
\begin{eqnarray*}
    \widetilde{\can}_A \circ \kappa(\suma a 1\sba \ot c\spa)
    &=& \suma a 1\sba \delta(c\spa{}\sw{1} \ot c\spa{}\sw{2}{}\sv{1}{}\zz{-1})
    c\spa{}\sw{2}{}\sv{1}{}\zz{0} c\spa{}\sw{2}{}\sv{2}{}\zz{0} \ot 
    c\spa{}\sw{2}{}\sv{2}{}\zz{1} \\
    &=& \sumab a 1\sba \delta(c\spa{}\sw{1} \ot c\spa{}\sw{2}) 1\sbb \ot 
    c\spa{}\sw{3}{}\spb \\
    &=& \sumab a 1\sba 1\sbb \ot c\spab = \suma a 1\sba \ot c\spa.
    \end{eqnarray*}
The second equality follows from equation \eqref{eq.coactns} in 
Lemma~\ref{lemma4}, while the third one is a consequence of
the definition of a cointegral. The final equality
follows from equation~\eqref{re1}.
   Thus $\kappa$ is a left $\cC$-colinear splitting of 
    $\widetilde{\can}_A$. So by Corollary~\ref{cor.mod.main}, $A$ is a
    principal $\cC$-comodule, hence $\can_A:A\ot_B A\to \cC$ is an isomorphism 
   of $A$-corings and
    $A$ is a projective left $B$-module. 
In order to show that $A$ is $C$-equivariantly projective as a left 
$B$-module we need to show that there exists a section
 of the left $B$-action 
   $_A\varrho:B \ot A \to A$ in $_B\mathcal{M}^C$, where 
   $_{A}\varrho$ is given by the
multiplication $\mu$ in $A$, i.e., 
a left $B$-module, right $C$-comodule map $\sigma$ 
such that $_A\varrho \circ \sigma = A$.
The projectivity of $A$ as a left $B$-module means that there exists a 
left $B$-linear map $\tilde\sigma: A \to B \ot A$ such that
$_A\varrho \circ \tilde\sigma = A$. Now we can construct 
\begin{equation} 
\sigma: A \to B \ot A, \qquad 
\sigma = (B \ot A \ot \delta) \circ (B \ot \varrho^A \ot C) \circ
(\tilde\sigma \ot C) \circ \varrho^A. 
\end{equation} 
The map $\sigma$ is a composition of left $B$-linear maps, therefore
left $B$-linear. We now show that $\sigma$ is a right $C$-colinear
section of $_A\varrho$. The fact  that $\sigma$ is right $C$-colinear
follows from the colinearity of the cointegral $\delta$, 
equation~\eqref{colin} (this is a standard argument for coseparable
coalgebras). Explicitly,
note that, writing $\tilde\sigma(a) = \sum a\su{1} \ot a\su{2}\in B\ot A$,
\begin{eqnarray*} 
    \varrho^{B \ot A} \circ \sigma(a) &=& \sum a\zz{0}\suc{1} \ot 
    a\zz{0}\suc{2}{}\zz{0}{}\zz{0} \ot a\zz{0}\suc{2}{}\zz{0}{}\zz{1} 
    \delta(a\zz{0}\suc{2}{}\zz{1} \ot a\zz{1})\\
&=&\sum a\zz{0}\suc{1} \ot 
    a\zz{0}\suc{2}{}\zz{0} \ot a\zz{0}\suc{2}{}\zz{1}{}\sw{1} 
    \delta(a\zz{0}\suc{2}{}\zz{1}{}\sw 2 \ot a\zz{1})\\
&=& \sum a\zz{0}{}\su{1} \ot a\zz{0}{}\su{2}{}\zz{0} 
	    \delta(a\zz{0}{}\su{2}{}\zz{1} \ot a\zz{1}{}\sw{1}) \ot 
	    a\zz{1}{}\sw{2} \\
 &=& \sum a\zz{0}{}\zz{0}\suc{1} \ot 
	a\zz{0}{}\zz{0}\suc{2}{}\zz{0} \delta(a\zz{0}{}\zz{0}\suc{2}{}\zz{1} 
	\ot a\zz{0}{}\zz{1}) \ot a\zz{1}\\
 &=&(\sigma \ot C) \circ \varrho^{A}(a).
    \end{eqnarray*}
	We can express the fact that $\tilde\sigma$ is a section of ${}_A\varrho$
	in the current notation as $\sum a\su{1} a\su{2} = a$. Using this property, 
 the fact that, by the definition of the coinvariants $B$, the right $C$-coaction
 is left $B$-linear,  
	and the fact that $\delta$ is a cointegral of $C$, we obtain
	\begin{eqnarray*}
	    {}_A\varrho(\sigma(a)) &=& \sum a\zz{0}\suc{1} 
    a\zz{0}\suc{2}{}\zz{0}  
    \delta(a\zz{0}\suc{2}{}\zz{1} \ot a\zz{1})\\
&=&	  \sum (a\zz{0}\suc{1} 
    a\zz{0}\suc{2}){}\zz{0}  
    \delta( (a\zz{0}\suc{1} 
    a\zz{0}\suc{2}){}\zz{1} \ot a\zz{1})\\
&=& \sum a\zz 0{}\zz 0\delta(a\zz 0{}\zz 1\ot a\zz{1}) = \sum  a\zz 0{}\delta(a\zz 1{}\sw 1\ot a\zz{1}\sw 2)\\
&=&  \sum a\zz 0\eps_C(a\zz 1) = a.
   \end{eqnarray*}
   This proves that $\sigma$ is a left $B$-module, right $C$-comodule section of the
   multiplication map ${}_A\varrho$ as required, and thus completes the proof of 
   the theorem.
\end{proof}
   
   In view of Corollary~\ref{cor.doi.inv}, Theorem~\ref{thm.main} implies
   the following
   \begin{corollary}\label{cor.main}
   Let  $H$ be a coseparable, $k$-projective weak Hopf algebra 
   with bijective antipode and let
   $A$ be a $k$-projective right $H$-comodule algebra. Let $\cE$ be the coring
   defined as the image of the projection $p_R$ given by 
\eqref{pr.Hopf}. If the map
   $$
\tcan_A :A \otimes A \to \cE, \qquad a
\otimes b \mapsto \sum ab \zz 0 \otimes b
\zz 1
$$
  is surjective, then  $A$ is a weak Hopf-Galois $H$-extension of the coinvariants 
  $B$ and  it is
$H$-equivariantly projective as a left $B$-module. 
\end{corollary}
  
  \section{Weak coalgebra Galois extensions with coalgebras projective as comodules}
  The aim of this section is to prove that, within an invertible weak entwining structure,
  the projectivity of $C$ as a $C$-comodule coupled with the surjectivity of the
  canonical map leads to a weak coalgebra-Galois extension. As a corollary we obtain a weak Hopf algebra version of the Kreimer-Takeuchi theorem 
  \cite[Theorem~1.7]{KreTak:hop} that states that for a right comodule algebra $A$ of a finite-dimensional Hopf algebra, the surjectivity of the canonical map implies that $A$ is a Hopf-Galois extension, and that $A$ is projective over its coinvariants.
  \begin{theorem}\label{thm.main.proj}
Let $(A,C,\re,\le)$ be an invertible weak entwining structure and let $\cC$ be a coring corresponding to the weak entwining $(A,C,\re)$ as in 
Proposition~\ref{prop.ent.cor}. Suppose that:
\begin{blist}
\item $A$ is a right weak entwined module with coaction $\rho^A$;
\item $C$ is $k$-flat and projective as a left $C$-comodule;
\item $(A \ot A)^{\mathrm{co}C} = A \ot B$, where $B=A^{\mathrm{co}C}$;
\item $\tcan_A: A \ot A \to \cC$, $a
\otimes b \mapsto \sum ab \zz 0 \otimes b
\zz 1$
 is surjective.
\end{blist}
Then $B\subseteq A$ is a weak coalgebra-Galois extension and $A$ is $k$-relatively projective as a left $B$-module.
\end{theorem}
\begin{proof}
First, notice that from Remark~\ref{rem.abelian} we can deduce that the assumption (c) above implies that $\lhom\cC A {A\ot A} = {}^{\mathrm{co}\cC}(A \ot A)_g = B \ot A$, where $g=\sum 1\zz 0
\ot 1\zz 1$. More precisely, 
$$
(A \ot A)^{\mathrm{co}C} = \{\sum_i a_i\ot b_i\in A\ot A \; |\; \sum_i a_i\ot gb_i = \sum_i a_i\ot b_ig\},
$$
hence the usual twist map $a\ot b\mapsto b\ot a$ gives rise to an 
isomorphism ${}^{\mathrm{co}\cC}(A \ot A)_g \simeq (A \ot A)^{\mathrm{co}C}$. Consequently, the assumption (c) is equivalent to the 
statement  that ${}^{\mathrm{co}\cC}(A \ot A)_g = B\ot A$, i.e.,  as explained in Remark~\ref{rem.abelian}, the condition (b) in
Theorem~\ref{thm.mod.main} is fulfilled.
Let $\cD$ be a coring corresponding to $(A,C,\le)$.
By Proposition~\ref{prop.inv.cor}, the coring $\cC$ is isomorphic to $\cD$, hence 
$\cC$ is a left $\cD$-comodule
via $(\le \otimes _A \cC) \circ \DC$. The correspondence between left $\cD$-comodules 
and left entwined modules then yields that $\cC$ is a left entwined
module, hence, in particular, a left $C$-comodule.  For $\suma 1\sba \otimes c\spa \in
\cC$, the left $C$-coaction comes out as
$$ ^{\cC}\varrho (\sum_\alpha 1\sba \otimes c \spa) = \sum_{\alpha,E}
c\spa \sw 1 \mathrm{} _E \otimes 1\sba \mathrm{} ^E \otimes c\spa
\sw 2 .
$$
Note that
\begin{eqnarray}
\sum_\alpha c\sw 1 \otimes 1\sba \otimes c\sw 2 \spa
    &=& \sum_\alpha c\sw 1 \otimes \eps (c\sw 2 \spa) 1\sba \otimes c\sw3 \nonumber 
    \\
    &=& \sum_E c\sw 1 \otimes \eps (c \sw 2 \mathrm{} _E)1^E \otimes c
\sw3 \nonumber\\
    &=& \sum _E c\sw 1 \mathrm{} _E \otimes 1^E \otimes c \sw 2. \label{left.coa}
\end{eqnarray}
The first equality comes from \eqref{re4}, the second is a consequence of the part (c) of 
Definition~\ref{def.inv} of an invertible weak entwining structure,  and the 
final equality is implied by \eqref{le4}.  In view of this we can compute
\begin{eqnarray*}
^\cC\varrho(\sum_\alpha 1\sba \otimes c\spa )
    &=& \sum_{\alpha, \beta, E} c\sw 1 \spb \mathrm{} _E \otimes
    {1\sbab} ^E \otimes c \sw 2 \spa  \\
    &=& \sum_{\alpha, E} c\sw 1 \mathrm{} _E \otimes 1^E 1\sba
    \otimes c \sw 2 \spa  \\
    &=& \sum _{E,F} c \sw 1 _{FE} \otimes 1^E 1^F \otimes c \sw 2 \\
    &=& \sum _E c \sw 1 \mathrm{} _E \otimes 1 ^E \otimes c \sw
    2.\\
\end{eqnarray*}
The first equality follows by the properties of right-right weak
entwining structures in particular (\ref{re3}), the second by part (b)
of Definition~\ref{def.inv}, the third by the previous computation
(\ref{left.coa}) and the final equality by (\ref{le1}).  This implies that the
map
\begin{equation}\label{map}
\ell:C \to \cC, \qquad c \mapsto \suma 1\sba \otimes c\spa 
\end{equation}
is a left $C$-comodule map.  As $\tcan_A$ is a left $\cC$-comodule
map (hence a $\cD$-comodule map by Proposition~\ref{prop.inv.cor}), 
it is a left $C$-comodule map.  Since $\tcan_A$ is surjective
and $C$ is projective as a left $C$-comodule, the map
$$
\xymatrix{
\Lhom C C {A \otimes A} \ar[rrr]^{\Lhom C C {\tcan_A}} &&& \Lhom C C \cC} 
$$
is surjective.  This implies that there exists
$\hat{f}\in\Lhom C C {A \otimes A}$ such that
\begin{equation}\label{split}
 \tcan_A \circ \hat{f} =\ell
 \end{equation}
with $\ell$ given by (\ref{map}). Any left $C$-comodule map
$\hat{f} : C \to A \otimes A$ satisfies assumption \eqref{eq.colin} of 
Lemma~\ref{lemma3}, hence there is a left $\cC$-comodule map
$$ f: \cC \to A \otimes A, \qquad f(\sum_ \alpha a 1\sba \otimes c
\spa) = \sum _ \alpha a 1 \sba \hat{f} (c\spa).$$ Take any $a\in
A$ and $c \in C$ and compute
\begin{eqnarray*}
\tcan _A \circ f(\sum _\alpha a1\sba \otimes c\spa)
&=& \tcan _A (\sum_\alpha a 1\sba \hat{f}(c\spa))
= \sum_\alpha a 1\sba \tcan _ A (\hat{f}(c\spa))\\
&=& \sum_\alpha a 1 \sba \ell (c\spa) 
= \sum_{\alpha, \beta} a 1\sba 1\sbb \otimes c \spab
= \sum_\alpha a 1\sba \otimes c\spa. 
\end{eqnarray*}
Second equality follows by left linearity of $\tcan_A$, the third
by (\ref{split}), the fourth by (\ref{map}) and the last by the property
\eqref{re1} of a right-right weak entwining map.  
This means that $f$ is a left $\cC$-colinear section of
$\tcan_A$.  With additional assumptions this means Theorem~\ref{thm.mod.main}
can be applied to obtain desired result.
\end{proof}

Theorem~\ref{thm.main.proj} leads to the following weak Hopf algebra version of
the Kreimer-Takeuchi theorem \cite[Theorem~1.7]{KreTak:hop}.
\begin{corollary}\label{kre.tak}
Let $k$ be a field and let $H$ be a finite dimensional weak Hopf algebra over $k$. 
Let $A$ be a right $H$-comodule algebra and $\cE$ be the $A$-coring associated to the corresponding right-right weak entwining map $\re$ given by \eqref{Hopf.ent.r}. 
 If 
 $$
\tcan_A :A \otimes A \to \cE, \qquad a
\otimes b \mapsto \sum ab \zz 0 \otimes b
\zz 1
$$ 
 is surjective, then $B \subseteq A$ is a weak Hopf-Galois extension and $A$ is  projective as a left (and right) $B$-module.
\end{corollary}
\begin{proof}
First, \cite[Theorem~2.10]{Boh:wea} implies that $H$ has a bijective antipode, hence
by Corollary~\ref{cor.doi.inv}, $(A,H,\re,\le)$ is an invertible weak entwining structure
with $\re$ given by \eqref{Hopf.ent.r} and $\le$ given by \eqref{Hopf.ent.l}. Second, as
the dual $H^*$ of a finite dimensional weak Hopf algebra is a weak Hopf algebra,
\cite[Theorem~3.11]{Boh:wea} implies that $H^*$ is a quasi-Frobenius algebra (i.e.,
it is self-injective). Now combination of the Faith-Walker theorem \cite[Theorem~24.12]{Fai:alg2} that asserts that every injective module over a quasi-Frobenius algebra is projective and \cite[Theorem~1.3]{GomNas:qua} that states that a quasi-co-Frobenius coalgebra is projective as a comodule,  implies
that $H$ is a projective left (and right) $H$-comodule 
(cf.\ \cite[Remark~1.5]{GomNas:qua}). Since $k$ is a field, the condition
(c) in Theorem~\ref{thm.main.proj} is automatically satisfied. Thus Theorem~\ref{thm.main.proj} yields the required assertion. The right $B$-projectivity of
$A$ follows from the left-right symmetry (see the discussion at the end of this section).
\end{proof}

Every weak Hopf algebra is a Hopf algebroid over $R=\Im \pl$  and a right $H$-comodule algebra $A$ is also a comodule algebra over this 
Hopf-algebroid (cf.\ \cite[36.9, 37.15]{BrzWis:cor} or see the original papers \cite{EtiNik:dyn}, \cite{Szl:fin}, \cite{Sch:wea} and \cite{BrzCae:Doi} for more details).   A right-right weak entwining structure $(A,H,\re)$ with $\re$ given by \eqref{Hopf.ent.r} can be understood as a right-right  entwining structure over $R$. With this identification, one can also deduce Corollary~\ref{kre.tak} from \cite[Corollary~4.3]{Boh:Gal}.

 By a theorem of Lin (cf.\  \cite[Proposition~5]{Lin:sem}),  every right co-Frobenius coalgebra is projective as a left comodule, hence  Theorem~\ref{thm.main.proj} can be applied to
co-Frobenius coalgebras (or, even more generally, to 
quasi-co-Frobenius coalgebras of G\'omez-Torrecillas and N\v ast\v asescu, cf.\ 
 \cite[Theorem~1.3]{GomNas:qua}). In particular, one obtains in this way a weak Hopf
algebra version of the Beattie-D\v asc\v alescu-Raianu extension of the Kreimer-Takeuchi theorem to co-Frobenius Hopf algebras 
(cf.\ \cite[Theorem~3.1]{BeaDas:Gal}). 

Finally let us mention that, throughout, we worked in a {\em right-handed convention},
taking right comodules (such as right Galois-comodules), right coalgebra extensions
(algebras $A$ with a right $C$-coaction), etc. Clearly, all the results presented
here can also be presented in a {\em left-handed convention}. For example, assuming
that the canonical map $\tcan_M$ in Theorem~\ref{thm.mod.main} is a split
epimorphism in the category of right $\cC$-comodules, we can deduce that $M^*$
is a left Galois $\cC$-comodule. In the case of invertible weak entwining structures
the distinction between left- and right-handed conventions is blurred in the sense that
every right weak coalgebra-Galois extension (corresponding to a right-right weak entwining $\re$) is a left weak coalgebra-Galois extension (corresponding to $\le$) and
vice versa. In particular this means that the assumption that $C$ is projective
as a left $C$-comodule in Theorem~\ref{thm.main.proj} can be replaced by the assumption that $C$ is a projective right $C$-comodule. As a consequence, in this case, one obtains
that $A$ is a weak coalgebra-Galois extension, $k$-relatively projective as a 
right $B$-module. The same arguments together with the left-right symmetry of the notions of a quasi-Frobenius algebra and a quasi-co-Frobenius coalgebra imply that
$A$ in Corollary~\ref{kre.tak} is projective as a right $B$-module.

\section*{Acknowledgements}
We thank Gabriella B\"ohm for comments. 
Tomasz Brzezi\'nski thanks the Engineering and Physical Sciences
Research
Council for an Advanced Fellowship.

\end{document}